\documentclass[11pt]{amsart}

\usepackage[official]{eurosym}
\usepackage{array,float,graphicx,epic,eepic}
\usepackage{amssymb, amsmath, amsthm, mathrsfs}

\theoremstyle{plain}
\newtheorem{thmx}{Theorem}[section]

\newtheorem{lemma}{Lemma}[section]
\newtheorem{theorem}{Theorem}[section]
\newtheorem{corollary}{Corollary}[section]
\newtheorem{proposition}{Proposition}[section]

\theoremstyle{definition}

\newtheorem{example}{Example}[section]
\newtheorem{remark}{Remark}[section]

\newenvironment{NotMyTheorem}{\bigskip \noindent{\bf Theorem}
\nonumber \it}{\bigskip}

\newenvironment{NotMyLemma}{\bigskip \noindent{\bf Lemma}
\nonumber \it}{\bigskip}

\newtheorem*{acknowledgements}{Acknowledgements}

\def \ol {\overline}
\def \Q {\ensuremath{\mathbb{Q}}}

\def \N {\ensuremath{\mathbb{N}}}

\def \mcM {\ensuremath{\mathcal{M}}}

\def \mcO {\ensuremath{\mathcal{O}}}
\def \mcR {\ensuremath{\mathcal{R}}}
\def \mcS {\ensuremath{\mathcal{S}}}
\def \mcV {\ensuremath{\mathcal{V}}}

\def \Qp {\mathbb{Q}_p}
\def \Z {\mathbb{Z}}

\def \Zp  {\mathbb{Z}_p}

\def \Fp  {\mathbb{F}_p}
\def \Fq  {\mathbb{F}_q}

\def \bfb {{\bf b}}

\def \bfD {{\bf D}}
\def \bff {{\bf f}}
\def \bfg {{\bf g}}

\def \bfk {{\bf k}}
\def \bfl {{\bf \ell}}
\def \bfm {{\bf m}}
\def \bfwtm {{\bf \widetilde{m}}}

\def \bfn {{\bf n}}
\def \mcN {\ensuremath{\mathcal{N}}}

\def \bfr {{\bf r}}
\def \bfs {{\bf s}}
\def \bfwts {{\bf \widetilde{s}}}
\def \bft {{\bf t}}
\def \bfX {{\bf X}}
\def \bfx {{\bf x}}
\def \bfY {{\bf Y}}
\def \bfy {{\bf y}}
\def \bfY {{\bf Y}}
\def \bfgamma {\boldsymbol{\gamma}}
\def \bfrho {\boldsymbol{\rho}}
\def \bfsigma {\boldsymbol{\sigma}}

\def \L {\Lambda}
\def \GL {\text{GL}}

\def \wt {\widetilde}
\def \u {\"{u}}

\def \mcI {\mathcal{I}}

\def \tud {\textup{d}}
\def \zetacc {\zeta^{\rm{cc}}}
\def \cc {\mathrm{cc}}
\def \LieL {\mathcal{L}}
\def \Adstar {{\rm{Ad}^*}}
\def \Ind {{\rm Ind}}
\def \Ltwos {{L'_s}}
\def \primcharndprime {(\Z/p^N)^{n}\setminus(p\Z/p^N)^{n}}
\def \T {\mathcal{T}}
\def \Ttwo {\mathcal{T}_2}
\def \zirr {\zeta^{\rm{irr}}}
\def \G {\Gamma}
\def \Rad {\textup{Rad}}
\def \tl {\triangleleft}
\def \matnzeroone {\text{Mat}_n(\{0,1\})}
\def \nl {\triangleleft}

\def \m {\phantom}

\begin{document}

\title[Functional equations for zeta functions of groups and
rings]{Functional equations for zeta functions\\ of groups and rings}

\date{\today}
\author{Christopher Voll}
\address{Email: C.Voll.98 at cantab.net }
\address{School of Mathematics, University of Southampton, Highfield, SO17 1BJ, United Kingdom.}


\begin{abstract} 
  
  We introduce a new method to compute explicit formulae for various
  zeta functions associated to groups and rings. The specific form of
  these formulae enables us to deduce local functional equations. More
  precisely, we prove local functional equations for the subring zeta
  functions associated to rings, the subgroup, conjugacy and
  representation zeta functions of finitely generated, torsion-free
  nilpotent (or $\T$-)groups, and the normal zeta functions of
  $\T$-groups of class~$2$. In particular we solve the two problems
  posed in~\cite[Section 5]{duSG/06}. We deduce our theorems from a
  `blueprint result' on certain $p$-adic integrals which generalises
  work of Denef and others on Igusa's local zeta function.  The Malcev
  correspondence and a Kirillov-type theory developed by Howe are used
  to `linearise' the problems of counting subgroups and
  representations in $\T$-groups, respectively.

\end{abstract}
\maketitle

\bigskip 

\section{Introduction}

Zeta functions of groups were introduced by Grunewald, Segal and Smith
in the 1980s as a tool to study the subgroup growth of finitely
generated groups. In \cite{GSS/88} the zeta function of a finitely
generated group $G$ was defined as the Dirichlet series
\begin{equation}\label{definition zeta}
\zeta_{G}(s)=\sum_{H\leq G} |G:H|^{-s},
\end{equation}
where $s$ is a complex variable and the sum ranges over the finite
index subgroups of $G$. Grunewald, Segal and Smith derived results on
the zeta functions of finitely generated, torsion-free nilpotent (or
$\T$-)groups~$G$, and went on to consider variants
of~\eqref{definition zeta}. These include a $\T$-group's \emph{normal}
zeta function counting only normal subgroups of finite index, and the
\emph{conjugacy} zeta function counting subgroups up to
conjugacy. They also developed an analogous theory for rings. (In the
current paper, by a ring we mean a finitely generated abelian group
with a bi-additive product.) The \emph{ideal} zeta function of a ring,
for instance, generalises the classical Dedekind zeta function of a
number field. Much of the subsequent developments in the theory of
zeta functions of groups and rings is documented in the
monograph~\cite{LubotzkySegal/03} and in the report~\cite{duSG/06}.

Only comparatively recently, Hrushovski and Martin
(\cite{MartinHrushovski/04}) started to investigate
\emph{representation} zeta functions of $\T$-groups, enumerating
(twist-isoclasses of) finite-dimensional complex characters.

All the zeta functions mentioned so far have the property that they
satisfy an \emph{Euler product} decomposition into local factors,
indexed by the primes. For example, for a $\T$-group $G$ we have
$$\zeta_{G}(s)=\prod_{p \text{ prime}}\zeta_{G,p}(s),$$ where
$\zeta_{G,p}(s)=\sum_{H\leq_p G}|G:H|^{-s}$ enumerates finite
$p$-power index subgroups. All these local zeta functions are known to
be rational functions in the parameter~$p^{-s}$ with integer
coefficients~(\cite{GSS/88},~\cite{MartinHrushovski/04}). In many
cases they exhibit a remarkable symmetry: For instance, it had been
observed (cf., e.g., \cite{duSG/06}) that the local factors of all
zeta functions of $\T$-groups $G$ for which explicit formulae are
known satisfy a \emph{local functional equation} of the form
\begin{equation}\label{funeq prelim}
\zeta_{G,p}(s)|_{p\rightarrow p^{-1}}=(-1)^a
p^{b-cs}\zeta_{G,p}(s),
\end{equation}
for almost all primes $p$ and suitable integers $a,b,c$ depending only
on the Hirsch length of $G$. Here $p\rightarrow p^{-1}$ denotes a
formal inversion of the local parameter $p$, which we shall now
explain. 

The prime example is the case of $G=\Z^n$. It is known from
\cite[Proposition 1.1]{GSS/88} that
$$\zeta_G(s) = \prod_{i=0}^{n-1}\zeta(s-i),$$ where
$\zeta(s)=\prod_{p \text{ prime}}\frac{1}{1-p^{-s}}$ is the Riemann
zeta function. The functional equation
$$\prod_{i=0}^{n-1}\frac{1}{1-p^{-(i-s)}}=(-1)^np^{\binom{n}{2}-ns}\prod_{i=0}^{n-1}\frac{1}{1-p^{i-s}}$$
of the local factor at the prime $p$ is easily seen to hold for all
primes. Here, as well as in all other cases in which explicit formulae
are known to date, the local zeta functions $\zeta_{G,p}(s)$ are in
fact rational functions in $p^{-s}$ \emph{and} $p$, and the left hand
side of~\eqref{funeq prelim} denotes the rational function obtained by
formally inverting both of these two parameters. This `uniformity' in
the prime~$p$, however, is not typical: Results of du Sautoy and
Grunewald~(\cite{duSG/00}) show that the dependence of the local
(normal) zeta functions of $\T$-groups on the primes will, in general,
reflect the variation of the number of $\Fp$-points of certain
algebraic varieties defined over~$\Fp$, which may be far from
polynomial in the prime~$p$ (see also~\cite{duS-ecI/01} and
\cite{duS-ecII/01}).  In~\cite{Voll/05} we produced examples of normal
zeta functions of $\T$-groups of nilpotency class~$2$ (or
$\Ttwo$-groups) which exhibit functional equations similar
to~\eqref{funeq prelim} and which are not uniform. Indeed, the local
factors of these zeta functions are rational functions in $p^{-s}$,
whose coefficients involve the numbers $b_V(p)$ of $\Fp$-rational
points of certain smooth projective varieties $V$ over $\Fp$ which
are, in general, not polynomials in~$p$. By the Weil conjectures,
these numbers may be expressed as alternating sums of Frobenius
eigenvalues. The operation $p\rightarrow p^{-1}$ is performed by
inverting these eigenvalues. In the special (`uniform') case that the
$b_V(p)$ are in fact polynomials in $p$, this specialises to an
inversion of the prime $p$.


Only a single example of a representation zeta function of a
$\T$-group seems to have appeared in print so far: In
\cite{MartinHrushovski/04} Hrushovski and Martin derive a formula for
the representation zeta function of the discrete Heisenberg group in
terms of the Riemann zeta function and its
inverse~(cf.~Example~\ref{example heisenberg}).

\medskip
According to du~Sautoy and Segal, to find an explanation for the
phenomenon of local functional equations for zeta functions of groups
and rings is ``one of the most intriguing open problems in this area''
(\cite[p.~274]{duSSegal/00}). In the current paper we prove that local
functional equations hold for (almost all factors of) the following
zeta functions:
\begin{itemize}
\item[(A)] zeta functions of rings (and, as a corollary,
of $\T$-groups),
\item[(B)] conjugacy zeta functions of $\T$-groups,
\item[(C)] normal zeta functions of $\Ttwo$-groups and
\item[(D)] representation zeta functions of $\T$-groups.
\end{itemize}

By proving (A) and (C) we solve Problems~5.1 and~5.2 posed
in~\cite{duSG/06}. In its given generality, (C) is best possible: It
is known that the normal zeta functions of nilpotent groups of
class~$3$ may or may not satisfy local functional
equations~(cf.~\cite{duSautoyWoodward/06}). To determine the exact
scope of this intriguing symmetry for ideal zeta functions of rings
remains a challenging open problem.
 
We achieve our results by showing that all of the above-mentioned zeta
functions may be expressed in terms of certain $p$-adic integrals,
generalising Igusa's local zeta function. Given a non-constant
polynomial $f(\bfy)\in\Z[y_1,\dots,y_m]$, its associated Igusa local
zeta function is the $p$-adic integral
$$\int_{\Zp^m}|f(\bfy)|^s|\tud\bfy|,$$ where $s$ is a complex
variable, $|\phantom{0}|$ stands for the $p$-adic absolute value and
$|\tud\bfy|$ denotes the (additive) Haar measure on $\Zp^m$, the
affine $m$-space over the $p$-adic integers~$\Zp$. This $p$-adic
integral is closely related to the Poincar\'e series counting $p$-adic
points on the hypersurface defined by~$f$ (cf.~\cite{Denef/91}). We
prove functional equations for these integrals by generalising results
by Denef and others on Igusa's local zeta function.

The integrals considered in the present paper are quite different from
the `cone integrals' introduced by du Sautoy and Grunewald
in~\cite{duSG/00} (see Section~\ref{section methodology} for further
details).

\subsection{Detailed statement of results}

Let $L$ be a ring. Its (subring) zeta function is defined to be the
Dirichlet series
$$\zeta_L(s) = \sum_{H\leq L} |L:H|^{-s},$$ where the sum ranges over
all subrings $H$ of finite index in $L$, and $s$ is a complex
variable. This zeta function decomposes naturally as an Euler product,
indexed by the primes:
$$\zeta_L(s) = \prod_{p \text{ prime}}\zeta_{L,p}(s),$$ where
$\zeta_{L,p}(s)=\zeta_{L\otimes\Zp}(s)$. Grunewald, Segal and Smith
proved in \cite{GSS/88} that each local factor is a rational function
in $p^{-s}$ with integral coefficients.  Our first main Theorem is

\begin{thmx}\label{theorem rings} 
Let $L$ be a ring of torsion-free rank $n$. Then there are smooth
projective varieties $V_t$, $t\in \{1,\dots,m\}$, defined over~$\Q$,
and rational functions $W_t(X,Y)\in\Q(X,Y)$ such that for almost all
primes~$p$ the following hold.
\begin{enumerate}
\item[(a)] Denoting by $b_t(p)$ the number of $\Fp$-rational points of
$\ol{V_t}$, the reduction$\mod p$ of~$V_t$, we have
\begin{equation}\label{formula rings}
\zeta_{L,p}(s) = \sum_{t=1}^m b_t(p) W_t(p,p^{-s}).
\end{equation}

\item[(b)] Setting $b_t(p^{-1}) := p^{-\dim(V_t)}b_t(p)$ the following
functional equation holds:\m{funeq rings}
\begin{equation}\label{funeq rings}
\zeta_{L,p}(s)|_{p\rightarrow p^{-1}} = (-1)^n
p^{\binom{n}{2}-ns}\zeta_{L,p}(s).
\end{equation}
\end{enumerate}
\end{thmx}

The novelty of this result is that it allows us to deduce the
equations~\eqref{funeq rings}. In \cite{duSG/00}, du Sautoy and
Grunewald gave formulae akin to~\eqref{formula rings} for, \emph{inter
alia}, the local factors of $\zeta_{L}(s)$. Their proof depends on a
representation of local zeta functions through certain $p$-adic
integrals called `cone integrals' which in general will not satisfy
functional equations like~\eqref{funeq rings}. See~\cite[Chapter
4]{duSautoyWoodward/06} for a discussion of functional equations for
cone integrals.
 
Given a ring $L$ we cannot, in general, pin down the primes $p$ which
have to be excluded in Theorem~\ref{theorem rings}. On the other hand,
any prime $p$ will be among the primes for which Theorem~\ref{theorem
rings}, applied to the ring $pL$, makes no assertion.

We also note that the definition in part (b) of Theorem~\ref{theorem
rings} is consistent with our above explanation of the operation
$p\rightarrow p^{-1}$. Indeed, by the Weil conjectures, the numbers
$b_t(p)$ may be expressed as alternating sums of Frobenius
eigenvalues. These complex numbers satisfy certain symmetries which
are reflected by the functional equations for the Weil zeta functions
of the varieties~$\ol{V_t}$. Had the expressions $b_t(p^{-1})$ been
defined as the numbers obtained from inverting these eigenvalues, the
identities $b_t(p^{-1}) = p^{-\dim(V_t)}b_t(p)$ would follow from
these symmetries (see the remarks preceding Theorem~\ref{inversion
property Z tilde} for details).


\begin{example}\label{example 3D}
Theorem~\ref{theorem rings} applies in particular to the zeta
functions of `simple' Lie algebras over $\Z$ such as
$\mathfrak{sl}(d,\Z)$. The only such Lie algebra for which a
functional equation as in~\eqref{funeq rings} had been previously
established is $L=\mathfrak{sl}(2,\Z)$
(cf.~\cite{duSTaylor/02}\footnote{The denominator of $P(2^{-s})$
in~\cite[Equation (7)]{duSG/06} should read $1-2^{1-3s}$.}):
$$\zeta_{L}(s)=\zeta(s)\zeta(s-1)\zeta(2s-1)\zeta(2s-2)(1+3\cdot2^{1-2s}-2^{3-3s})\prod_{p\not=2}(1-p^{1-3s}).$$
In fact, $\mathfrak{sl}(2,\Z)$ seems to be the only non-soluble Lie
ring whose subring zeta function has been computed explicitly. Note
that the functional equation fails for $p=2$. A similar phenomenon may
occur if one studies the zeta function of a $\Zp$-algebra:
In~\cite{Klopsch/03} Klopsch computes the zeta function of a maximal
$\Zp$-order in a central simple $\Qp$-division algebra of
index~$2$. The fact that this zeta function does not satisfy a
functional equation of the form~\eqref{funeq rings} reflects the fact
that it is not the `generic' local factor of the zeta function of a
ring.

Klopsch and the present author have unified and generalised the above
examples. In~\cite{KlopschVoll/07} they gave a formula for the zeta
function of an arbitrary $3$-dimensional $\Zp$-Lie algebra, based on
the proof of Theorem~\ref{theorem rings}.
\end{example}

An important corollary of Theorem~\ref{theorem rings} is to the theory
of zeta functions of finitely generated, torsion-free nilpotent (or
$\T$-) groups. In \cite[Theorem 4.1]{GSS/88} it was shown that, given
a $\T$-group $G$ of Hirsch length $n$, there is a Lie subring
$L=L(G)$, lying as a full $\Z$-lattice in an $n$-dimensional Lie
algebra $\LieL(G)$ over~$\Q$ such that for almost all
primes~$p$\m{group ring}
\begin{equation}\label{group ring}\zeta_{G,p}(s)=\zeta_{L,p}(s).
\end{equation}
 Thus we obtain

\m{corollary subgroups}
\begin{corollary}\label{corollary subgroups} 
  Let $G$ be a $\T$-group of Hirsch length $n$. For all but finitely
  many primes $p$
  $$\zeta_{G,p}(s)|_{p\rightarrow p^{-1}} = (-1)^n
  p^{\binom{n}{2}-ns}\zeta_{G,p}(s).$$
\end{corollary}

Theorem~\ref{theorem rings} is itself an instance of an application of
a `reciprocity' result (Corollary~\ref{corollary 2 to inversion
property} to Theorem~\ref{inversion property Z tilde}) establishing
certain functional equations for a family of $p$-adic integrals.
Theorem~\ref{inversion property Z tilde} may be viewed as a
generalisation of Stanley's `reciprocity theorem for linear
homogeneous diophantine equations'~(\cite[Theorem
4.6.14]{Stanley/98}), which we now briefly explain. Given a set of
simultaneous linear homogeneous diophantine equations in $n$
indeterminates, say, one may encode their non-negative (positive)
solutions in a rational generating function $E(\bfx)$ ($\ol{E}(\bfx)$,
respectively), where $\bfx=(x_1,\dots,x_n)$ is a vector of formal
variables. More precisely, one defines
$$E(\bfx):=\sum_{\alpha\in\N_0^n\cap\mathscr{C}}\bfx^\alpha \text{ and
}\ol{E}(\bfx):=\sum_{\alpha\in\N^n\cap\mathscr{C}}\bfx^\alpha,$$ where
$\alpha=(\alpha_1,\dots,\alpha_n)$, $\bfx^\alpha=x_1^{\alpha_1}\dots
x_n^{\alpha_n}$ and $\mathscr{C}$ is the cone of all non-negative
\emph{real} solutions to the given set of equations. Stanley's
reciprocity theorem states that, if $\ol{E}(\bfx)\not=0$, then
$$\ol{E}(1/\bfx)=(-1)^{\dim{\mathscr{C}}}E(\bfx).$$

Our generalisation is obtained using (a variant of) Stanley's result
and an explicit formula of the form~\eqref{formula rings} for the
integrals in question (Corollary~\ref{corollary denef projective} of
Theorem~\ref{theorem denef projective}). This formula in turn is
inspired by and generalises work of Denef (\cite{Denef/87}), Denef and
Meuser (\cite{DenefMeuser/91}) and Veys and Zuniga-Galindo
(\cite{VeysZG/06}) on Igusa's local zeta function. While Igusa's local
zeta functions associated to homogeneous polynomial mappings may be
expressed as integrals over projective space, the $p$-adic integrals
considered in the current paper reduce to integrals over the complete
flag variety $\text{GL}_n/B$, where $B$ is a Borel subgroup. In a
sense this explains the factor $p^{\binom{n}{2}}$ in~\eqref{funeq
rings}.

\bigskip 
A variant of the problem of counting subgroups consists in counting
subgroups only up to conjugacy. Let $G$ be a $\T$-group. The conjugacy
zeta function of a $\T$-group $G$ is defined as \m{cc definition}
\begin{equation*}\label{cc definition} 
\zetacc_{G}(s)=\sum_{H\leq G} |G:H|^{-s}|\mathcal{C}_G(H)|^{-1}
\end{equation*}
where $|\mathcal{C}_{G}(H)|$ is the size of the conjugacy class
of~$H$. It is known (\cite[Remark on p.~189]{GSS/88}) that
$\zetacc_G(s)$ also has an Euler product decomposition into local
factors $\zetacc_{G,p}(s)$ which are all rational in~$p^{-s}$. By
applying the results of Section~\ref{section igusa} we shall prove
\m{theorem conjugacy}
\begin{thmx}
\label{theorem conjugacy}
Let $G$ be a $\T$-group of Hirsch length $n$. For all but finitely
many primes~$p$
$$\zetacc_{G,p}(s)|_{p\rightarrow p^{-1}}=(-1)^n
p^{\binom{n}{2}-ns}\zetacc_{G,p}(s).$$
\end{thmx}

\bigskip As a third application of our rather technical `blueprint
result' Theorem~\ref{theorem denef projective} and its applications we
deduce functional equations for normal zeta functions of
$\Ttwo$-groups. The normal zeta function of a $\T$-group $G$ is
defined as
$$\zeta^\triangleleft_{G}(s)=\sum_{H\triangleleft G}|G:H|^{-s}$$ where
the sum ranges over the \emph{normal} subgroups $H$ of finite index in
$G$. It also satisfies an Euler product decomposition. We prove

\begin{thmx}\label{theorem normal} 
  Let $G$ be a $\Ttwo$-group of Hirsch length $n$ with centre $Z(G)$
  such that $G/Z(G)$ has torsion-free rank~$d$. For all but finitely
  many primes $p$
  $$\zeta^\triangleleft_{G,p}(s)|_{p\rightarrow p^{-1}} = (-1)^n
  p^{\binom{n}{2}-(d+n)s}\zeta^\triangleleft_{G,p}(s).$$
\end{thmx}
As mentioned above, a functional equation may or may not hold for
normal zeta functions of class greater than
three. See~\cite[Chapter~2]{duSautoyWoodward/06} for examples and
\cite[Theorem 4.44]{duSautoyWoodward/06} for a conjectural form in
case it does hold.

\bigskip The fourth and last application in this paper of
Theorem~\ref{theorem denef projective} and its consequences is
concerned with \emph{representation} zeta functions of $\T$-groups,
which we shall now explain. Given a $\T$-group $G$, we denote by
$R_n(G)$ the set of $n$-dimensional complex characters of~$G$. Given
$\sigma_1,\sigma_2\in R_n(G)$, we say that $\sigma_1$ and $\sigma_2$
are {\sl twist--equivalent} if there exists a linear
character~$\chi\in R_1(G)$ such that~$\sigma_1=\chi\sigma_2$. The
classes of this equivalence relation are called twist-isoclasses. We
say that a character $\sigma$ of a representation $\rho$ of $G$
factors through a finite quotient of $G$ if $\rho$ factors through it.
The set $R_n(G)$ has the structure of a quasi--affine complex
algebraic variety whose geometry was analysed by Lubotzky and Magid
in~\cite{LubotzkyMagid/85}. They proved

\begin{NotMyTheorem}\cite[Theorem 6.6]{LubotzkyMagid/85} \label{theorem
    luma} Let $G$ be a $\T$-group. For every~$n\in\N$ there is a
  finite quotient $G(n)$ of $G$ such that every $n$-dimensional
  irreducible character of~$G$ is twist-equivalent to one that
  factors through~$G(n)$. In particular, the number of
  twist-isoclasses of irreducible $n$-dimensional characters is
  finite.
\end{NotMyTheorem}

Let us call this number $a_n$. The representation zeta function of $G$
is defined (cf.~\cite{MartinHrushovski/04}) by
$$\zirr_{G}(s):=\sum_{n=1}^\infty a_nn^{-s}.$$ It follows from the
above theorem and~\cite[(10.33)]{CurtisReiner-methods/81} that the
function~$n\mapsto a_n$ is multiplicative and thus
$$\zirr_{G}(s)=\prod_{p \text{ prime}}\zirr_{G,p}(s),$$
where
$$\zirr_{G,p}(s):=\sum_{n=0}^\infty a_{p^n}p^{-sn}.$$ 

\m{example heisenberg}
\begin{example}\label{example heisenberg} (\cite[Example 8.12]{MartinHrushovski/04}, \cite[Theorem 5]{NunleyMagid/89}) Let
$$H=\langle x_1,x_2,y| [x_1,x_2]=y, \text{ all other $[\,,]$
trivial}\rangle$$ be the discrete Heisenberg group. Then
$$\zirr_H(s)=\sum_{n=1}^\infty\phi(n)n^{-s}=\frac{\zeta(s-1)}{\zeta(s)}=\prod_{p\text{
prime}}\frac{1-p^{-s}}{1-p^{1-s}},$$ where $\phi$ denotes the Euler
totient function.
\end{example}

By a model-theoretic result of Hrushovski and
Martin~\cite[Theorem~8.4]{MartinHrushovski/04}, the local
representation zeta functions of a $\T$-group are known to be rational
functions in $p^{-s}$ with integer coefficients. By expressing
$\zirr_{G,p}(s)$ in terms of $p$-adic integrals to which
Theorem~\ref{theorem denef projective} is applicable we shall prove

\begin{thmx}\label{theorem representations} Let $G$ be a $\T$-group
  with derived group $G'=[G,G]$ of Hirsch length $n$. Then for almost
  all primes $p$
\begin{equation}\label{funeq representations}
\nonumber\zirr_{G,p}(s)|_{p\rightarrow p^{-1}} =  p^n\zirr_{G,p}(s).
\end{equation}
\end{thmx}

\subsection{Outline of methodology and related work}\label{section methodology}

We briefly describe how we relate the problems solved in
Theorems~\ref{theorem rings}, \ref{theorem conjugacy}, \ref{theorem
normal} and \ref{theorem representations} to problems about $p$-adic
integrals.

To count the subrings of $p$-power index in a ring~$L$ of rank~$n$, we
observe that it is enough to keep track of the index of the largest
$\Zp$-subalgebra of $L\otimes\Zp$ in each given homothety class of
lattices in the $p$-adic vector space~$\Qp^n$. Using the action of the
group $\G:=\text{GL}_n(\Zp)$ on the set of homothety classes, we show
that the latter problem reduces to counting polynomial congruences in
finite quotients of $\G$. This counting problem translates into the
problem of computing a $p$-adic integral in very much the same fashion
as the problem of counting polynomial congruences in affine space
translates to the problem of computing Igusa's local zeta function
(cf.~\cite[Section 1.2]{Denef/91}). It proved helpful to think of
homothety classes of lattices as the vertices of the Bruhat-Tits
building of $\text{SL}_n(\Qp)$, and to partition the vertex set into
finitely many parts according to their position relative to the `root
class' $[L\otimes\Zp]$. As we remarked above, this approach differs
decisively from the `cone integrals' introduced by du Sautoy and
Grunewald in~\cite{duSG/00}. Their analysis rests on a basis dependent
parametrisation of $p$-power index subrings of a given ring in terms
of upper-triangular matrices over the $p$-adic integers satisfying
certain divisibility conditions (`cone conditions').

To count subgroups up to conjugacy in a $\T$-group~$G$ we use the fact
that, for almost all primes $p$,
$$\zetacc_{G,p}(s)=\zetacc_{L,p}(s),$$ where $L=L(G)$ is the Lie ring
associated to $G$ (cf. the remark preceding Corollary~\ref{corollary
subgroups}), and
$$\zetacc_{L,p}(s)=\sum_{H\leq
L\otimes\Zp}|L\otimes\Zp:H|^{-s}|L\otimes\Zp:\mathcal{N}_{L\otimes\Zp}(H)|^{-1},$$
where $H$ ranges over the subalgebras of $L\otimes\Zp$of finite index
and $\mathcal{N}_{L\otimes\Zp}(H)$ is the normaliser of $H$
in~$L\otimes\Zp$. We thus have to keep track both of the largest
subring of $L\otimes\Zp$ in each given homothety class and of the
class' normaliser. The index of the latter is given by the index of a
system of linear congruences. Enumerating these indices, in turn, may
be achieved by counting the elementary divisors of matrices of linear
forms, encoding the group's commutator structure. In
Proposition~\ref{proposition elementary divisors} we show that,
slightly more generally, the generating functions enumerating
elementary divisors of matrices of polynomial forms of the same degree
may be expressed in terms of $p$-adic integrals associated to
degeneracy loci of these matrices, to which Corollary~\ref{corollary 1
to inversion property} is applicable.

In order to count normal subgroups in class-$2$-nilpotent groups we
develop an idea first introduced in~\cite{Voll/04}. There it was shown
that it suffices to evaluate a weight function on the set of homothety
classes of lattices in the centre $Z(L\otimes\Zp)$ of the $\Zp$-Lie
algebra $L\otimes\Zp$, where $L$ is the associated Lie ring. The
weight associated to a vertex in the appropriate affine Bruhat-Tits
building corresponding to a given class is again given by the index of
a system of linear congruences.

As mentioned above, the validity of our main Theorems~\ref{theorem
rings} and~\ref{theorem normal} had been known in many special cases.
We refer the reader to the numerous examples collected
in~\cite{duSautoyWoodward/06}. In this research monograph du~Sautoy
and Woodward also present a conjecture on functional equations for
cone integrals that would explain functional equations for normal zeta
functions of nilpotent of class greater than~$2$.

The key to Theorem~\ref{theorem representations} is to use
`Kirillov-theory' developed by Howe~\cite{Howe-nilpotent/77} to
translate the problem of counting irreducible representations of~$G$
to the problem of counting co-adjoint orbits in the dual of the Lie
algebra associated to $G$ by the Malcev correspondence. We use the
fact -- also established by Howe -- that the sizes of co-adjoint
orbits may be expressed in terms of the indices of the radicals of
certain anti-symmetric forms on the Lie algebra. These may also be
described in terms of elementary divisors of matrices encoding the
structure of the Lie algebra.

Representation zeta functions of nilpotent groups have not been
studied until fairly recently, and~\cite{MartinHrushovski/04} seems to
be the only reference so far on this topic.  The idea of using
Kirillov-theory to study representation zeta functions of groups,
however, has been successfully employed before. Jaikin-Zapirain proved
in~\cite{Zapirain/06} the rationality of representation zeta functions
for certain compact $p$-adic analytic groups using a Kirillov-type
correspondence developed by Howe~(\cite{Howe-compact/77}) for these
groups.  In~\cite{MartinHrushovski/04} Hrushovski and Martin suggest
that Jaikin-Zapirain's work may be adapted to prove rationality of
local representation zeta functions for $\T$-groups, too.

\medskip
Among the variants of the zeta function~\eqref{definition zeta} of a
$\T$-group $G$ considered in~\cite{GSS/88} is also the zeta function
$\zeta^{\wedge}_G(s)$, enumerating subgroups of finite index whose
profinite completion is isomorphic to the profinite completion
of~$G$. In~\cite{duSLubotzky/96} du~Sautoy and Lubotzky proved a
functional equations for the local factors of $\zeta^{\wedge}_G(s)$
for a class of $\T$-groups. Their work is based on a reduction of the
problem of computing these zeta functions to the problem of computing
certain $p$-adic integrals over the group's algebraic automorphism
group, generalising work of Igusa's~(\cite{Igusa/89}). The functional
equation for the local factors of these zeta functions arises from a
symmetry in the root system of the associated Weyl groups. An argument
of this kind (albeit only for the Weyl groups of type~$A$) is also
used in the present paper to deduce Corollary~\ref{corollary 2 to
inversion property}. We do not know whether the zeta functions
$\zeta^{\wedge{}}_G(s)$ (for reasonably large classes of $\T$-groups)
may be described by the $p$-adic integrals studied in the current
paper. In \cite{Berman/06} Berman extends the approach taken
in~\cite{duSLubotzky/96}, proving uniformity and local functional
equations for these zeta functions for a wider class of nilpotent
groups than previously considered.

\medskip
The fact that we have to disregard finitely many primes in most of our
results has two reasons: Firstly, our Theorems~\ref{theorem denef
affine} and~\ref{theorem denef projective} upon which
Theorem~\ref{inversion property Z tilde} and its corollaries are based
are valid only for primes for which a certain principalisation of
ideals has good reduction. Secondly, we are forced to ignore finitely
many primes in order to transfer between the $\T$-group $G$ and its
associated Lie algebra.

\subsection{Layout of the paper and notation}

In Section ~\ref{section igusa} we first develop, in
  Theorems~\ref{theorem denef affine} and~\ref{theorem denef
  projective}, explicit formulae for certain families of $p$-adic
  integrals generalising Igusa's local zeta function.  These two
  results -- which may be understood as close analogues of Theorems~2
  and~3 in~\cite{DenefMeuser/91} -- form the technical core of the
  paper. We use them to establish, in Theorem~\ref{inversion property
  Z tilde}, an `inversion property' enjoyed by the $p$-adic integrals
  considered. In Corollaries~\ref{corollary 2 to inversion property}
  and~\ref{corollary 1 to inversion property} we exploit this property
  to deduce functional equations for certain linear combinations of
  the $p$-adic integrals in question. The remainder of the paper is
  dedicated to showing how Theorem~\ref{theorem denef projective} may
  be used as a template to describe various kinds of zeta functions.
  In the second part of Section~\ref{section igusa} we give a first
  application of this idea to the problem of counting elementary
  divisors of matrices of forms (Proposition~\ref{proposition
  elementary divisors}).  In the four subsections of
  Section~\ref{section applications} we prove Theorems~\ref{theorem
  rings}, \ref{theorem conjugacy}, \ref{theorem normal} and
  \ref{theorem representations}, respectively.

 We use the following notation.
\medskip

\begin{tabular}{l|l}
  $\N$ & the set $\{1,2,\dots\}$ of natural numbers \\
$I=\{i_1,\dots,i_l\}_<$ & the set $I$ of natural numbers
$i_1<\dots<i_l$\\ $I_0$ & the set $I\cup\{0\}$ for $I\subseteq\N$ \\
$[k]$ & the set $\{1,\dots,k\}$, $k\in\N$ \\ $[l,k]$ & the set
$\{l,\dots,k\}$, $k,l\in\N$ \\ 
$\binom{a}{b}$ & the binomial coefficient for $a, b \in
\N_0$\\ $\binom{a}{b}_X$ & the polynomial $\prod_{i = 0}^{b-1} (1 -
X^{a-i})/(1 - X^{b-i})$, \\ & \quad where $a,b \in \N_0$ with $a \geq
b$ \\ & Note: The \emph{$q$-binomial coefficient} or
\emph{Gaussian}\\&\quad\emph{polynomial} $\binom{a}{b}_q$ gives the
number of \\ & \quad subspaces of dimension $b$ in $\Fq^a$.\\
$\binom{n}{I}_X$ & the polynomial
$\binom{n}{i_l}_X\binom{i_l}{i_{l-1}}_X\dots\binom{i_2}{i_1}_X$,\\ &
\quad for $n \in \N$, $I = \{i_1,\dots,i_l\}_< \subseteq[n-1]$\\ &
Note: $\binom{n}{I}_q\text{ gives the number of flags of type~$I$
in~$\Fq^n$.}$\\ 
$S_n$ & the symmetric
group on $n$ letters\\ $M^t$& the transpose of a matrix $M$\\
$\Zp$ & the ring of $p$-adic integers ($p$ a prime)\\
$\Qp$ & the field of $p$-adic numbers\\
$[\Lambda]$& the homothety class $\Qp^*\Lambda$ of a (full) lattice
$\Lambda$ in $\Qp^n$\\
$K$ & a finite extension of the field $\Qp$\\
$R$ & the valuation ring of $K$\\
$P$ & the maximal ideal of $R$\\
$\ol{K}$ & the residue field $R/P$, of cardinality $q$\\
$F$& a number field\\
$\delta_P$ & the `Kronecker delta' which is equal to $1$ if\\ & the property $P$ holds and equal to $0$ otherwise.
\end{tabular}
\medskip

\noindent Given a set $\bff$ of polynomials and a polynomial $g$, we
write $g\bff$ for $\{gf|\,f\in\bff\}$, and $(\bff)$ for the polynomial ideal generated by~$\bff$.

\section{Functional equations for some $p$-adic integrals}\label{section igusa}
\subsection{A blueprint result}
In this section we study a family of $p$-adic integrals generalising
Igusa's local zeta functions. We first introduce some more
notation. Let $p$ be a prime and $K$ be a finite extension of the
field $\Qp$ of $p$-adic numbers. Let $R=R_K$ denote the valuation ring
of $K$, $P=P_K$ the maximal ideal of $R$, and $\ol{K}$ the residue
field $R/P$. The cardinality of $\ol{K}$ will be denoted by $q$.

For $x\in K$, let $v(x)=v_{P}(x)\in\Z\cup\{\infty\}$ denote the
$P$-adic valuation of~$x$, and $|x|:=q^{-v(x)}$. For a finite set
$\mcS$ of elements of $K$, we set
$\|\mcS\|:=\max\{|s|\,|\;s\in\mcS\}$. Fix $k,m,n\in\N$. For each
$\kappa\in[k]$, let $(\bff_{\kappa\iota})_{\iota\in I_\kappa}$ be a
finite family of finite sets of polynomials in $K[y_1,\dots,y_m]$, and
let $x_1,\dots,x_{n-1}$ be independent variables. Also, for
$i\in[n-1]$ we fix non-negative integers $e_{i\,\kappa\iota}$. For a
set $I=\{i_1,\dots,i_l\}_<\subseteq[n-1]$, $\kappa\in[k]$, we set
\begin{equation*}
\bfg_{\kappa,I}(\bfx,\bfy)=\bigcup_{\iota\in I_\kappa}\left(\prod_{i\in I}x_i^{e_{i\,\kappa\iota}}\right)\bff_{\kappa\iota}(\bfy).
\end{equation*}
Let $W\subseteq R^m$ be a subset which is a union of cosets$\mod P^m$
and $\bfs=(s_1,\dots,s_k)$ be independent complex variables. We then
define
\begin{equation}\label{definition igusa}
  Z_{W,K,I}(\bfs):=\int_{P^l\times
  W}\prod_{\kappa\in[k]}\|\bfg_{\kappa,I}(\bfx,\bfy)\|^{s_\kappa}|\tud\bfx_I||\tud\bfy|
\end{equation}
where $|\tud\bfx_I|=|\tud x_{i_1}\wedge\dots\wedge\tud x_{i_l}|$ is
the Haar measure on $K^l$ normalised so that $R^l$ has measure~$1$
(and thus $P^l$ has measure~$q^{-l}$), and $|\tud\bfy|=|\tud
y_1\wedge\dots\wedge\tud y_m|$ is the (normalised) Haar measure on
$K^m$. It is well-known that $Z_{W,K,I}(\bfs)$ is a rational function
in $q^{-s_\kappa}$, $\kappa\in[k]$, with integral coefficients.

We now assume that the polynomials constituting the sets
$\bff_{\kappa\iota}$ are in fact defined over a number field $F$. We
may consider the local zeta functions $Z_{W,K,I}(\bfs)$ for all
non-archimedean completions $K$ of $F$. In the remainder of this
section we shall derive formulae for $Z_{W,K,I}(\bfs)$, valid for
almost all completions $K$ of $F$ under this and further
assumptions. They are essentially based on the formulae Denef gave for
Igusa's local zeta function
$$Z(s)=\int_{R^m}|f(\bfy)|^s|\tud\bfy|$$
in~\cite[Theorem~3.1]{Denef/87}, using the concept of resolution of
singularities for the hypersurface defined by $f$. In the case where
the single polynomial $f$ is replaced by a finite set of polynomials
$\bff$, Veys and Zuniga-Galindo (\cite[Theorem 2.10]{VeysZG/06}) gave
an analogous formula, using instead the concept of principalisation of
ideals, which we briefly recall.

\begin{NotMyTheorem}\cite[Theorem 1.0.1]{Wlodarczyk/05}\label{theorem principalisation}
  Let $\mcI$ be a sheaf of ideals on a smooth algebraic
  variety~$X$. There exists a principalisation $(Y,h)$ of $\mcI$, that
  is, a sequence
  $$
  X=X_0\stackrel{h_1}{\longleftarrow}X_1\longleftarrow\dots\stackrel{h_\iota}{\longleftarrow}X_{\iota}\longleftarrow\dots\stackrel{h_r}\longleftarrow
  X_r=Y$$ of blow-ups $h_{\iota}:X_{{\iota}}\rightarrow X_{\iota-1}$
  of smooth centres $C_{{\iota}-1}\subset X_{{\iota}-1}$ such that
\begin{itemize}
\item[a)] The exceptional divisor $E_{\iota}$ of the induced morphism
  $h^{\iota}=h_{\iota}\circ\dots\circ h_{1}:X_{\iota}\rightarrow X$
  has only simple normal crossings and $C_{\iota}$ has simple normal
  crossings with~$E_{\iota}$.
\item[b)] Setting $h=h_r\circ\dots\circ h_1$, the total transform
  $h^{*}(\mcI)$ is the ideal of a simple normal crossing
  divisor~$\wt{E}$. If the subscheme determined by $\mcI$ has no
  components of codimension one, then $\wt{E}$ is an $\N$-linear
  combination of the irreducible components of the divisor~$E_r$.
\end{itemize}
\end{NotMyTheorem}

Also recall the definition \cite[Definition 2.2 (\emph{mutatis
mutandis})]{Denef/87} of a \emph{principalisation $(Y,h)$ with good
reduction$\mod P$} if $\mcI$ and $(Y,h)$ are defined over a $p$-adic
field $K$. Note that, given a principalisation $(Y,h)$ for $\mcI$
defined over a number field $F$, $(Y,h)$ will have good reduction$\mod
P$ (where $P$ is the maximal ideal in the ring of integers of the
completion of $F$ at $K$) for almost all maximal ideals~$P$ of the
ring of integers of $F$ (this is essentially \cite[Theorem
2.4]{Denef/87}).

Specifically, let $(Y,h)$, $h:Y\rightarrow \mathbb{A}^m$, be a
principalisation of the ideal
$$\mcI=\prod_{\kappa\in[k],\,\iota\in I_\kappa}(\bff_{\kappa\iota})$$
where $(\bff)$ denotes the ideal generated by the finite set $\bff$ of
polynomials.  We set $\mcV:=\rm{Spec}(F[\bfy]/\mcI)$ and
$\mcV_{\kappa\iota}:=\rm{Spec}(F[\bfy]/(\bff_{\kappa\iota}))$. Then,
denoting by $E_t$, $t\in T$, the irreducible components of
$(h^{-1}(\mcV))_{\rm{red}}$, we have
\begin{align}
h^{-1}(\mcV)=&\sum_{t\in T}N_t E_t,\label{def N iota}\\
 h^{-1}(\mcV_{\kappa\iota})=&\sum_{t\in T}N_{t\,\kappa\iota}
 E_t,\label{def N iota bff}
\end{align}
say, for suitable non-negative integers $N_t,
N_{t\,\kappa\iota}$. Note that, for every $t\in T$,
$$N_t = \sum_{\kappa\in[k],\,\iota\in
I_\kappa}N_{t\,\kappa\iota}.$$ Similarly we denote by $\nu_t-1$ the
multiplicity of $E_t$ in the divisor of $h^*(\tud y_1\wedge\dots\wedge
\tud y_m)$. The numbers $(N_{t\,\kappa\iota},\nu_t)_{\substack{t\in
T,\,\kappa\in[k],\,\iota\in I_\kappa}}$ will be called the \emph{numerical
data} of the principalisation $(Y,h)$.\m{theorem denef affine}

\begin{theorem} \label{theorem denef affine} 
Suppose that all the sets $\bff_{\kappa\iota}$ are integral
(i.e. contained in $R[\bfy]$) and do not define the zero ideal$\mod
P_K$, and that $(Y,h)$ has good reduction$\mod P_K$. Then
\begin{equation*}
Z_{W,K,I}(\bfs)=\frac{(1-q^{-1})^{|I|}}{q^{m}}\sum_{U\subseteq
T}c_{U,W}(q)(q-1)^{|U|}\Xi_{U,I}(q,\bfs),
\end{equation*}
where
\begin{equation*}
c_{U,W}(q)=|\{a\in\ol{Y}(\ol{K})|\;a\in\ol{E_u}(\ol{K})\Leftrightarrow u\in U \text{ and }\ol{h}(a)\in \ol{W}\}|
\end{equation*}
(where $\overline{\phantom{0}}$ denotes reduction$\mod P$ and
$\ol{W}=\{(\ol{y_1},\dots,\ol{y_m})|\;(y_1,\dots,y_m)\in W\}|$) and
\begin{equation}\label{def Z U}
\Xi_{U,I} (q,\bfs)=\sum_{\substack{(m_u)_{u\in
U}\in\N^{|U|}\\(n_i)_{i\in I}\in\N^l}}q^{-\sum_{i}n_i-\sum_{u}\nu_u
m_u-\sum_{\kappa}s_\kappa\min_{\iota\in
I_\kappa}\{\sum_{i}e_{i\,\kappa\iota}n_i+\sum_{u}N_{u\,\kappa\iota}m_u\}}.\end{equation}
\end{theorem}

\begin{example} 
If $l=0$ and, for each $\kappa\in[k]$, $|I_\kappa|=1$,
Theorem~\ref{theorem denef affine} reduces to (a multivariable version
of) Vey's and Zuniga-Galindo's generalisation
\cite[Theorem~2.10]{VeysZG/06} to polynomial mappings of Denef's
explicit formula \cite[Theorem~2]{Denef/87} for Igusa's local zeta
function associated to a single polynomial. Notice in particular that
in this case 
$$\Xi_{U,\varnothing}(q,\bfs)=\sum_{(m_u)_{u\in
U}\in\N^{|U|}}q^{\sum_um_u(-\nu_u-\sum_{\kappa}N_{u\kappa}s_\kappa)}=\prod_{u\in
U}\frac{X_u}{1-X_u}$$ for
$X_u:=q^{-\nu_u-\sum_{\kappa}N_{u\kappa}s_\kappa}$, where we write $N_{u\kappa}$
for $N_{u\,\kappa\iota}$. Also compare Example~\ref{example 3} for the other `extremal
case' $\Xi_{\varnothing,I}(q,\bfs)$.
\end{example}

\begin{proof}
  The proof is analogous to the one of \cite[Theorem 3.1]{Denef/87}
  (=\cite[Theorem 2]{DenefMeuser/91}), with the concept of resolution
  of singularities replaced by the concept of principalisation of
  ideals. We adopt -- \emph{mutatis mutandis} -- Denef's notation and
  just explain how the proof differs from his. Let $a$ be a closed
  point of $\ol{Y}$, and thus also of $\wt{Y}$. Let $T_a=\{t\in
  T|\;a\in\ol{E_t}\}=\{t_1,\dots,t_r\}_<$, say. Define $H=\{b\in
  Y(K)|\;h(b)\in R^m\}$ and recall the definition of the
  `reduction$\mod P$'-map $\theta:H\rightarrow \ol{Y}(\ol{K})$. In the
  regular local ring $\mcO_{\wt{Y},a}$, there exist irreducible
  elements $\gamma_1,\dots,\gamma_m$ such that, on $\theta^{-1}(a)$, for
  all $\kappa\in[k]$, $\iota\in I_\kappa$,
\begin{align*}
\|\bff_{\kappa\iota}\circ h\| =&
|\gamma_1|^{N_{t_1\,\kappa\iota}}\dots
|\gamma_r|^{N_{t_r\,\kappa\iota}}\text{ and }\\ |h^*(\tud
y_1\wedge\dots\tud
y_m)|=&|\gamma_1|^{\nu_{t_1}-1}|\dots|\gamma_r|^{\nu_{t_r}-1}|\tud\gamma_1\wedge\dots\wedge\tud\gamma_m|.
\end{align*}
Setting $\tud\bfgamma:=\tud\gamma_1\wedge\dots\wedge\tud\gamma_m$ we
define
\begin{align*}
\lefteqn{Z_{a,I}(\bfs)}\\
:=&\int_{P^l\times\theta^{-1}(a)}\prod_{\kappa\in[k]}\max_{\iota\in
I_\kappa}\left\{\prod_{i\in
I}|x_i|^{e_{i\,\kappa\iota}}\prod_{\rho\in[r]}|\gamma_\rho|^{N_{t_\rho\,\kappa\iota}}\right\}^{s_\kappa}\prod_{\rho\in[r]}|\gamma_\rho|^{\nu_{t_\rho}-1}|\tud\bfx_I||\tud\bfgamma|\\
=&\int_{P^{l+m}}\prod_{\kappa\in[k]}\max_{\iota\in
I_\kappa}\left\{\prod_{i\in
I}|x_i|^{e_{i\,\kappa\iota}}\prod_{\rho\in[r]}|y_\rho|^{N_{t_\rho\,\kappa\iota}}\right\}^{s_\kappa}\prod_{\rho\in[r]}|y_\rho|^{\nu_{t_\rho}-1}|\tud\bfx_I||\tud\bfy|\\
=&\frac{(q-1)^{r+l}}{q^{m+l}}\sum_{\substack{(m_t)_{t\in
T_a}\in\N^{r}\\(n_i)_{i\in I}\in\N^l}}q^{-\sum_{i}n_i-\sum_{t}\nu_t
m_t-\sum_{\kappa}s_\kappa\min_{\iota\in
I_\kappa}\{\sum_{i}e_{i\,\kappa\iota}n_i+\sum_{t}N_{t\,\kappa\iota}m_t\}}.
\end{align*}
This suffices as $Z_{W,K,I}(\bfs)=
\sum_{\substack{a\in\ol{Y}(\ol{K})\\\ol{h}(a)\in
\ol{W}}}Z_{a,I}(\bfs)$.
\end{proof}

We now make the further assumption that $m=n^2$. We identify $K^{n^2}$
with $\text{Mat}_n(K)$ and assume that the ideals
$(\bff_{\kappa\iota})$, $\kappa\in[k]$, $\iota\in I_\kappa$, are
$B(F)$-invariant, where $B(F)$ is the group of $F$-rational points of
the Borel subgroup of upper-triangular matrices in $G=\text{GL}_n$,
acting on $K[y_{11}, y_{12},\dots,y_{nn}]$ by matrix-multiplication
from the right. Let $(Y,h)$, $h:Y\rightarrow G/B$ be a
principalisation of the ideal
$\mcI=\prod_{\kappa,\iota}(\bff_{\kappa\iota})$.  Denoting, similar to
the above, by $\mcV$ the subvariety of $G/B(K)$ defined by $\mcI$ and
by $\mcV_{\kappa\iota}$ the subvariety defined by
$(\bff_{\kappa\iota})$ yields numerical data
$(N_{t\,\kappa\iota},\nu_t)_{\substack{t\in T,\,\kappa\in[k],\,\iota\in
I_\kappa}}$ defined by formulae analogous to \eqref{def N iota} and
\eqref{def N iota bff} above. We shall study the integral
$$Z_I(\bfs):=Z_{W,K,I}(\bfs)$$ for $W=\G=\GL_n(R)$ for almost all
completions~$K$ of~$F$. Note that the Haar measure $\mu'$ on the
compact topological group $\G$ coincides with the additive Haar
measure $\mu$ induced from $R^{n^2}$ (and normalised such that
$\mu(R^{n^2})=1$), as $\mu'=|\det|^{-n}\mu=\mu$. This will be
important in later applications as it implies, for example, that all
the cosets of a finite index subgroup $\G'\leq\G$ have measure
$\mu(\G)/|\G:\G'|$, with
$\mu(\G)=(1-q^{-1})\dots(1-q^{-n})$. \m{theorem denef projective}

\begin{theorem}\label{theorem denef projective} Suppose that, in addition to the above assumptions, none of the ideals $(\bff_{\kappa\iota})$ is equal to the zero ideal$\mod P_K$, and that $(Y,h)$ has good reduction$\mod P_K$. Then
\begin{equation*}\label{formula denef projective} 
Z_I(\bfs)=\frac{(1-q^{-1})^{|I|+n}}{q^{\binom{n}{2}}}\sum_{U\subseteq
T}c_{U}(q)(q-1)^{|U|}\Xi_{U,I}(q,\bfs),
\end{equation*}
where each $c_U(q)$ is the number of $\ol{K}$-rational points of
$\ol{E_U}\setminus\cup_{V\supsetneq U}\ol{E_V}$ ($E_U:=\cap_{u\in
U}E_u$) and $\Xi_{U,I}(q,\bfs)$ is defined as in \eqref{def Z U}
above.
\end{theorem}

\begin{proof}
  The proof follows closely the spirit of the proof of \cite[Theorem
  3]{DenefMeuser/91}. In fact, our function $Z_I(\bfs)$ is a close
  analogue of the function $\widehat{Z}_K(s)$, defined on~\cite[p.
  1140]{DenefMeuser/91}. We write $\G$ as a disjoint union of sets
  $$\G_\sigma=\{\bfx\in \G|\;\ol{\bfx}\in B(\Fq)\sigma B(\Fq)\},$$
  $\sigma\in S_n$, where $\text{GL}_n(\Fq)=\bigcup_{\sigma\in
  S_n}B(\Fq)\sigma B(\Fq)$ is the Bruhat decomposition (here
  $\sigma\in S_n$ is identified with the respective permutation matrix
  in $\text{GL}_n(\Fq)$).  Thus $$Z_I(\bfs)=\sum_{\sigma\in
  S_n}Z_{\G_\sigma,K,I}(\bfs).$$ There is an obvious map
  $\gamma:\G\rightarrow G/B(K)$, and, by our invariance assumption on
  the ideals $(\bff_{\kappa\iota})$, the value of the integrand of
  $Z_I(\bfs)$ at a point $(\bfx,\bfy)\in P^l\times \G$ only depends on
  $\bfx$ and $\gamma(\bfy)$. By taking the measure $\omega$ on
  $G/B(K)$ which induces the Haar measure on the unit ball
  $R^{\binom{n}{2}}$ of each affine chart satisfying $\omega(a +
  P^{\binom{n}{2}})=q^{-\binom{n}{2}}$ and noting that
  $\mu(B)=(1-q^{-1})^n$, we obtain
  $$Z_{\G_\sigma,K,I}(\bfs)=(1-q^{-1})^n\int_{P^l\times
    V_\sigma}\prod_{\kappa\in[k]}\|\bfg_{\kappa,I}(\bfx,\bfy)\|^{s_\kappa}|\tud\bfx_I|\tud\omega,$$
  where $V_\sigma=\gamma(\G_\sigma)$. The projective variety $G/B$ may
  be covered by varieties $U_\sigma$, isomorphic to affine
  $\binom{n}{2}$-space, indexed by the elements of the symmetric group
  $S_n$, such that each $V_\sigma$ is contained in $U_\sigma$ and is a
  union of cosets$\mod P^{\binom{n}{2}}$. Theorem~\ref{theorem denef
    affine} may thus be applied to the restriction
  $(Y^\sigma,h^\sigma)$ of $(Y,h)$, a principalisation of the ideal
  defining the restriction of $\mcV$ to $U_\sigma$
  ($Y^\sigma=h^{-1}(U^\sigma)$, $h^\sigma=h|_{Y^\sigma}$), with good
  reduction$\mod P$. We obtain
  $$\int_{P^l\times
    V_\sigma}\prod_{\kappa\in[k]}\|\bfg_{\kappa,I}(\bfx,\bfy)\|^{s_\kappa}|\tud\bfx_I|\tud\omega=\frac{(1-q^{-1})^l}{q^{\binom{n}{2}}}\sum_{U\subseteq
    T}c_{U,\sigma}(q)(q-1)^{|U|}\Xi_{U,I}(q,\bfs),$$
  where

  $$c_{U,\sigma}(q)=|\{a\in\ol{Y}^\sigma(\ol{K})|\;a\in\ol{E_u}(\ol{K})\Leftrightarrow
  u\in U\text{ and }\ol{h}(a)\in\ol{V_\sigma}\}|.$$ The result follows
  since, if $a\in\ol{Y}(\ol{K})$, then $\ol{h}(a)$ is in exactly one
  $\ol{V_\sigma}$. Thus $\sum_{\sigma\in
  S_n}c_{U,\sigma}(q)=|\{a\in\ol{Y}(\ol{K})|\;a\in\ol{E_u}(\ol{K})\Leftrightarrow
  u\in U\}|=c_U(q)$.
\end{proof}

We now consider the normalised integrals \m{definition Z I tilde}
\begin{equation}\label{definition Z I tilde}
\wt{Z_I}(\bfs):=\frac{Z_I(\bfs)}{(1-q^{-1})^{|I|}\mu(\G)}.
\end{equation}

\begin{corollary}\label{corollary denef projective} 
For $U\subseteq T$, let $b_U(q)$ denote the number of
  $\ol{K}$-rational points of $\ol{E_U}$. Then 
  \begin{equation}\label{formula Z tilde}
\wt{Z_I}(\bfs) = |G/B(\Fq)|^{-1}\sum_{U\subseteq
    T}b_U(q)\sum_{V\subseteq U}(-1)^{|U\setminus
    V|}(q-1)^{|V|}\Xi_{V,I}(q,\bfs).
\end{equation}
\end{corollary}

\begin{proof} 
This follows immediately from the formula given for $Z_I(\bfs)$ in
Theorem~\ref{theorem denef projective}, Definition~\eqref{definition Z
I tilde}, the fact that $|G/B(\Fq)|=\binom{n}{[n-1]}_q$ and from the
identity
  $$c_V(q)=\sum_{V\subseteq U\subseteq T}(-1)^{|U\setminus V|}b_U(q).$$
\end{proof}

Before we proceed we consider a very special case.\m{example 3}
\begin{example}\label{example 3}
Assume that, for all $\kappa\in[k]$, $|I_\kappa|=1$ and that all
$\bff_{\kappa\iota}=\{1\}$. We write $e_{i\kappa}$ for
$e_{i\,\kappa\iota}$. We have
$$\Xi_{\varnothing,I}(q,\bfs)=\sum_{(n_i)_{i\in I}\in\N^l}q^{\sum_i
n_i(-1-\sum_\kappa{e_{i\kappa}}s_\kappa)}=\prod_{i\in
I}\frac{X_i}{1-X_i},$$ with
$X_i=q^{-1-\sum_\kappa{e_{i\kappa}}s_\kappa}$. Also note that
$b_{\varnothing}(q)=|G/B(\Fq)|$ and that $b_U(q)=0$ if
$U\not=\varnothing$. Thus $\wt{Z_I}(\bfs)=\prod_{i\in
I}\frac{X_i}{1-X_i}$. It is trivial to verify the `inversion property'
\begin{equation}
\label{funeq example}\wt{Z_I}(\bfs)|_{q\rightarrow q^{-1}}=(-1)^{|I|}\sum_{J\subseteq I}\wt{Z_J}(\bfs),
\end{equation}
as
$$\prod_{i\in I}\frac{X_i^{-1}}{1-X_i^{-1}}=(-1)^{|I|}\sum_{J\subseteq
I}\prod_{j\in J}\frac{X_j}{1-X_j}.$$
\end{example}
In the remainder of the current section we shall show that
equation~\eqref{funeq example} holds under the premises of
Theorem~\ref{theorem denef projective}. To give meaning to the
left-most term in \eqref{funeq example} in general we have to explain
what we mean by $b_U(q^{-1})$ (the other constituents of the
expression \eqref{formula Z tilde} for $\wt{Z_I}(\bfs)$ being rational
functions in $q$ and $q^{-s_1},\dots,q^{-s_k}$). Recall that by
properties of the Weil zeta functions associated to the
$\left(\binom{n}{2}-|U|\right)$-dimensional smooth projective
varieties $\ol{E_U}$ it is known that
$$
b_U(q)=\sum_{r=0}^{2(\binom{n}{2}-|U|)}(-1)^r\sum_{j=1}^{t_{U,r}}\alpha_{U,r,j}
$$ for suitable non-negative integers~$t_{U,r}$ and non-zero complex
numbers~$\alpha_{U,r,j}$, with the property that, for each $U$,
$r$, the \underline{multi}sets
$$
\left\{\left.{\alpha_{U,r,j}}\right|j\in[t_{U,2(\binom{n}{2}-|U|)-r}]\right\}
$$ 
and
$$
\left\{\left.\frac{q^{\binom{n}{2}-|U|}}{\alpha_{U,r,j}}\right|j\in[t_{U,r}]\right\}
$$ coincide (cf., e.g., \cite[Proof of Theorem
4]{DenefMeuser/91}). This motivates the definition
\begin{equation}\label{b q inverse}
  b_U(q^{-1}):=q^{-\left(\binom{n}{2}-|U|\right)}b_U(q)=\sum_{r=0}^{2(\binom{n}{2}-|U|)}(-1)^r\sum_{j=1}^{t_{U,r}}\alpha^{-1}_{U,r,j}.
\end{equation}
 We shall prove

\begin{theorem} \label{inversion property Z tilde}\m{inversion property Z tilde}Under the assumptions of Theorem~\ref{theorem denef projective}, the
following `inversion properties' hold:
\begin{equation}\label{inversion property}\tag{IP} 
\forall I\subseteq[n-1]:\;\wt{Z_I}(\bfs)|_{q\rightarrow q^{-1}} =
(-1)^{|I|}\sum_{J\subseteq I}\wt{Z_J}(\bfs)
\end{equation}
\end{theorem}

\begin{proof}
To see what happens to the (rational) functions $\Xi_{V,I}(q,\bfs)$ in
expression~\eqref{formula Z tilde} if we formally invert the prime
power $q$, we employ a result of Stanley:

\begin{proposition}\label{proposition stanley}
  Let $L_{\sigma\tau}(\bfn)$, $\sigma\in[s]$, $\tau\in[t]$, be
  $\Z$-linear forms in the variables $n_1,\dots,n_r$ and
  $X_1,\dots,X_r,Y_1,\dots,Y_s$ independent variables, and set
\begin{align*}
  Z^\circ(\bfX,\bfY):=&\sum_{\bfn\in\N^r}\prod_{\rho\in[r]}
  X_\rho^{n_\rho}\prod_{\sigma\in[s]}Y_\sigma^{\min_{\tau\in[t]}\{L_{\sigma\tau}(\bfn)\}}\\
  Z(\bfX,\bfY):=&\sum_{\bfn\in\N_0^r}\prod_{\rho\in[r]}
  X_\rho^{n_\rho}\prod_{\sigma\in[s]}Y_\sigma^{\min_{\tau\in[t]}\{L_{\sigma\tau}(\bfn)\}}
\end{align*}
Then
$$Z^\circ(\bfX^{-1},\bfY^{-1}) = (-1)^r Z(\bfX,\bfY).$$
\end{proposition}
\begin{proof}
The proof of \cite[Theorem 4.6.14]{Stanley/98} carries through to this
slightly more general situation, provided one chooses a
triangulation of $\N^r$ that refines a subdivision into rational
polyhedral cones eliminating the `$\min$'-terms in the sum defining
$Z(\bfX,\bfY)$.
\end{proof}

\begin{corollary} 
  For all $I\subseteq[n-1]$, $V\subseteq T$,
\begin{equation}\label{corollary}
  \Xi_{V,I}(q,\bfs)|_{q\rightarrow q^{-1}} = (-1)^{|V|+|I|}\sum_{\substack{W\subseteq V, J\subseteq I}} \Xi_{W,J}(q,\bfs).
\end{equation}
\end{corollary}
We record the following simple fact:
\begin{lemma}
For all $U\subseteq T, J\subseteq[n-1]$,
\begin{multline}\label{umformung}
  \sum_{V\subseteq U}(-1)^{|U\setminus V|} (1-q^{-1})^{|V|}\sum_{W\subseteq V} \Xi_{W,J}(q,\bfs)=\\q^{-|U|}\sum_{V\subseteq U}(-1)^{|U\setminus V|}(q-1)^{|V|}\Xi_{V,J}(q,\bfs).
\end{multline}
\end{lemma}
The proof is a simple computation. We can now deduce
\begin{align*}
 \wt{Z_I}(\bfs)|_{q\rightarrow q^{-1}}=&\frac{(-1)^{|I|}q^{\binom{n}{2}}}{|G/B(\Fq)|}\sum_{U\subseteq
  T}b_U(q^{-1})\sum_{V\subseteq U}(-1)^{|U\setminus
  V|}(1-q^{-1})^{|V|}\cdot\\&\hspace{0.5in}\sum_{W\subseteq V, J\subseteq
  I}\Xi_{W,J}(q,\bfs)&&\eqref{formula Z tilde},\eqref{corollary}\\
 =&(-1)^{|I|}\sum_{J\subseteq
  I}|G/B(\Fq)|^{-1}\sum_{U\subseteq T}q^{|U|}b_U(q)\cdot\\&\hspace{0.5in}\sum_{V\subseteq
  U}(-1)^{|U\setminus V|}(1-q^{-1})^{|V|}\sum_{W\subseteq
  V}\Xi_{W,J}(q,\bfs)&&\eqref{b q inverse}\\
  =&(-1)^{|I|}\sum_{J\subseteq
  I}|G/B(\Fq)|^{-1}\sum_{U\subseteq T}b_U(q)\cdot\\&\hspace{0.5in}\sum_{V\subseteq
  U}(-1)^{|U\setminus V|}(q-1)^{|V|}\Xi_{V,J}(q,\bfs)&&\eqref{umformung}\\
  =&(-1)^{|I|}\sum_{J\subseteq
  I}\wt{Z_J}(\bfs).&&\eqref{formula Z tilde}
\end{align*}
This completes the proof of Theorem~\ref{inversion property Z tilde}.
\end{proof}

Recall that the polynomials $\binom{n}{I}_X$ were introduced at the
end of the introduction. We define
\begin{equation}\label{definition Z tilde}
\wt{Z}(\bfs) =
\sum_{I\subseteq[n-1]}\binom{n}{I}_{q^{-1}}\wt{Z_I}(\bfs).
\end{equation}
 
To prove Theorems~\ref{theorem rings},~\ref{theorem conjugacy}
and~\ref{theorem normal} we shall need\m{corollary 2 to inversion
property}
\begin{corollary}\label{corollary 2 to inversion property}
Under the assumptions of Theorem~\ref{theorem denef projective}, the following functional equation holds:
\begin{equation}\label{corollary 2}
\wt{Z}(\bfs)|_{q\rightarrow q^{-1}} = (-1)^{n-1}q^{\binom{n}{2}} \wt{Z}(\bfs).
\end{equation}
\end{corollary}

\begin{proof}This follows from the proof of~\cite[Corollary 2]{VollBLMS/06}. Note that Theorem~\ref{theorem denef projective} provides the required analogue of~\cite[Lemma~6]{VollBLMS/06}.
\end{proof}

Theorem~\ref{theorem representations} will follow from
Proposition~\ref{proposition elementary divisors} of the next section
which is in turn a special case of the following straightforward
corollary.

\begin{corollary}\label{corollary 1 to inversion property}
Under the assumptions of Theorem~\ref{theorem denef projective}, for
any $i\in[n-1]$,
\begin{equation}\label{corollary 1}
\left(\wt{Z_{\varnothing}}(\bfs) + (1-q^{-n})\wt{Z_{\{i\}}}(\bfs)\right)|_{q\rightarrow q^{-1}}=q^n \left(\wt{Z_{\varnothing}}(\bfs) + (1-q^{-n})\wt{Z_{\{i\}}}(\bfs)\right).
\end{equation}
\end{corollary}


\subsection{A first application: Counting elementary divisors}\label{subsection elementary divisors} We show how the problem of counting elementary divisors of matrices of
forms may be reduced to the problem of computing $p$-adic integrals of
the form studied in the previous section, associated to the polynomials
describing the degeneracy loci of these matrices. The main result of this subsection -- Proposition~\ref{proposition elementary divisors} -- will be needed to prove Theorem~\ref{theorem representations} in Section~\ref{subsection howe}.

Again, let $K$ be a $p$-adic field with valuation ring $R$, whose
maximal ideal is denoted by~$P$. Let $\mcR$ be an $e\times f$-matrix (with $e\geq f$, say) of
polynomials $\mcR_{ij}(\bfY)\in R[Y_1,\dots,Y_n]$. We make the
assumption on $\mcR$ that, whenever $\bfy=(y_1,\dots,y_n)\in R^n$ is a
vector with $\ol{\bfy}\not=0$, at least one entry of $\ol{\mcR(\bfy)}$
is non-zero. For a non-negative integer $N$ and a vector $\bfy
P^N\in(R/P^N)^n$ we say that $\mcR(\bfy P^N)$ has \emph{elementary
  divisor type} $\bfm$ (written $\nu(\mcR(\bfy P^N))=\bfm$) if
$\bfm=(m_1,\dots,m_f)$, $m_i\in[0,N]$, $m_1\leq\dots\leq m_f$, and
there are matrices $\beta\in\text{GL}_e(R/P^N)$,
$\gamma\in\text{GL}_f(R/P^N)$ such that
\begin{equation*}
  \beta\mcR(\bfy P^N)\gamma\equiv\left(\begin{array}{ccc}q^{m_1}&&\\&\ddots&\\&&q^{m_f}\\&&\end{array}\right)\mod P^N.
\end{equation*}
For $\bfm\in\N_0^f$ we set
$$ \mcN_{N,\mcR,\bfm}:=\left|\left\{ \bfy P^N \in (R/P^N)^n|\;\ol{\bfy
      P^N}\not=0,\nu(\mcR(\bfy P^N))=\bfm\right\}\right|.
$$ Note that $\mcN_{N,\mcR,\bfm}=0$ unless $0=m_1\leq\dots\leq m_f\leq
N$ (the necessity of $m_1=0$ being a consequence of our assumption
on~$\mcR$). Given, in addition, a $g\times h$-submatrix $\mcS$ of
$\mcR$ (WLOG $g\geq h$) defined by choosing $g$ rows and $h$ columns
of $\mcR$, and an $h$-tuple $\bfn$ we define
$$\mcN_{N,\mcR,\mcS,\bfm,\bfn}:=\left|\left\{ \bfy P^N \in(R/P^N)^n
|\;\ol{\bfy}\not=0,\begin{array}{c}\nu(\mcR(\bfy P^N))=\bfm,\\
  \nu(\mcS(\bfy P^N))=\bfn\end{array}\right\}\right|.$$ Again,
  $\mcN_{N,\mcR,\mcS,\bfm,\bfn}=0$ unless $0=m_1\leq\dots\leq m_f\leq
  N$ and $n_1\leq\dots\leq n_h$. We suppress the subscripts $\mcR$ and
  $\mcS$ if they are clear from the context.

Given complex variables $r_1,\dots,r_{f},s_1,\dots,s_{h}$, we define
the generating function
$$P(\bfr,\bfs)=P_{\mcR,\mcS,K}(\bfr,\bfs) =
\sum_{\substack{N\in\N_0\\\bfm\in\N_{0}^f,\bfn\in\N_{0}^h}}\mcN_{N,\bfm,\bfn}q^{-\sum_{i\in[f]}(N-m_i)r_i-\sum_{j\in[h]}(N-n_j)s_j}.$$

We now assume that the matrix $\mcR$ is in fact defined over a number
field~$F$, that its entries are all homogeneous of the same degree and
that the above assumption on $\mcR$ is satisfied for almost all
completions $K$ of $F$ for which all $\mcR_{ij}(\bfY)\in
R_K[\bfY]$. We consider such a `good' completion $K$ and drop the
subscript $K$. For $i\in[f]_0$, let $\bfrho_i$ denote the set of
$i$-minors of $\mcR$. The polynomials $\bfrho_i$ define the
$(\text{rk}\leq i-1)$-locus (or $i$-th degeneracy locus) of
$\mcR(\bfY)$. Similarly, let $\bfsigma_j$, $j\in[h]_0$, denote the set
of $j$-minors of $\mcS$.  Let
\begin{equation*}
  k:=\max\{i\in[f]_0|\;(\bfrho_i)\not=(0)\}\quad\text{ and }\quad
  l:=\max\{j\in[h]_0|\;(\bfsigma_j)\not=(0)\}.
\end{equation*}
Note firstly that $\bfrho_0=\bfsigma_0=\{1\}$, secondly that, by our
assumption on $\mcR$, $k\geq 1$, thirdly that $0\leq l\leq k$ and
fourthly that $P(\bfr,\bfs)$ is really a function in the variables
$r_1,\dots,r_k,s_1,\dots,s_l$:
\begin{equation}\label{P rep growth}
P(\bfr,\bfs)=\sum_{\substack{N\in\N_0\\\bfm\in\N_0^{k},
    \bfn\in\N_0^{l}}}\mcN_{N,\bfm,\bfn}q^{-\sum_{\kappa\in[k]}(N-m_\kappa)r_\kappa-\sum_{\lambda\in[l]}(N-n_\lambda)s_\lambda}
\end{equation}
where we set, given $\bfm=(m_1,\dots,m_k)$ and $\bfn=(n_1,\dots,n_l)$,
$$\mcN_{N,\bfm,\bfn}:=\mcN_{N,(m_1,\dots,m_k,N,\dots,N),(n_1,\dots,n_l,N,\dots,N)}.$$

For $I\subseteq\{1\}$ and $W=\G=\text{GL}_n(R)$ as above, consider the
$p$-adic integral
\begin{multline} 
Z_{I}(\bfr,\wt{\bfr},\bfs,\wt{\bfs},t):=\int_{P^{|I|}\times
\G}|x^{|I|}|^t\prod_{\kappa\in[k]}\|\bfrho_\kappa(\bfy^1)\cup
x^{|I|}\bfrho_{\kappa-1}(\bfy^1)\|^{r_\kappa}\|\bfrho_{\kappa-1}(\bfy^1)\|^{\wt{r}_\kappa}\cdot\\
\quad\prod_{\lambda\in[l]}\|\bfsigma_\lambda(\bfy^1)\cup
x^{|I|}\bfsigma_{\lambda-1}(\bfy^1)\|^{s_\lambda}\|\bfsigma_{\lambda-1}(\bfy^1)\|^{\wt{s}_\lambda}|\tud
x_I||\tud\bfy|,
\end{multline}
where $\bfy^1$ denotes the first column of the matrix $\bfy\in \G$.

\begin{remark}\label{remark elementary divisors} Whilst artificial, the formulation of $Z_I
$ as an integral over
$P^{|I|}\times \G$ rather than over $P^{|I|}\times R^n\setminus P^n$
serves to make it fit the `blueprint' Theorem~\ref{theorem denef projective} provided in the previous section. 
\end{remark}

 We now set, for $\bfm=(m_1,\dots,m_k)\in\N_0^k$,
 $\bfn=(n_1,\dots,n_l)\in\N_0^l$ and $N\in\N$
$$\mu_{N,\bfm,\bfn}=\mu\left\{(x,\bfy)\in P\times \G|\,v(x)=N,
\nu(\mcR(\bfy^1 P^N))=\bfm, \nu(\mcS(\bfy^1 P^N))=\bfn\right\} $$ and
$Z_I(\bfr,\bfs,t):=Z_I(\bfr,-\bfr,\bfs,-\bfs,t)$. Note that, by
definition, $\mu_{N,\bfm,\bfn}=0$ unless $0\leq m_1\leq\dots\leq
m_k\leq N$ and $n_1\leq\dots\leq n_l\leq N$. By definition of the
polynomials $\bfrho_{\kappa}$, $\bfsigma_{\lambda}$, we have
\begin{align}
Z_{\varnothing}(\bfr,\bfs,t)&=\mu(\G)\quad (\text{and thus,
  by~\eqref{definition Z I tilde},
  }\wt{Z_\varnothing}(\bfr,\bfs,t)=1),\label{Z empty rep growth}\\
  Z_{\{1\}}(\bfr,\bfs,t)&=\sum_{\substack{N\in\N\\\bfm\in\N_0^k,\bfn\in\N_0^l}}\mu_{N,\bfm,\bfn}q^{-tN-\sum_{\kappa}r_\kappa
  m_\kappa - \sum_{\lambda}s_\lambda n_{\lambda}}. \label{Z 1 rep
  growth}
\end{align}

Theorem~\ref{theorem denef projective} is applicable to $Z_I$,
$I\subseteq\{1\}$, together with a principalisation $(Y,h)$,
$h:Y\rightarrow G/B$, of the ideal
$\mcI=\prod_{\kappa\in[k]}(\bfrho_\kappa)\prod_{\lambda\in[l]}(\bfsigma_\lambda)$. Indeed,
$B(F)$-invariance is a consequence of the fact that the entries of
$\mcR$ were all assumed to be homogeneous of the same degree.

The following crucial lemma relates the numbers $\mu_{N,\bfm,\bfn}$
with the data $\mcN_{N,\bfm,\bfn}$ we would like to capture:\m{lemma
relation eldiv}
\begin{lemma}\label{lemma relation eldiv} For $N\in\N$
\begin{equation}\label{relation N mu rep growth}
\mcN_{N,\bfm,\bfn}=
\frac{\binom{n}{1}_{q^{-1}}}{\mu(\G)}\mu_{N,\bfm,\bfn}q^{N(n+1)}.
\end{equation}
\end{lemma}
Recall that $\binom{n}{1}_{q^{-1}}=(1-q^{-n})/(1-q^{-1})$.
\begin{proof}
Let $\G_{\{1\},N}$ denote the group 
$$\left\{\left(\begin{array}{c|c}\gamma_1&*\\\hline
q^N*&\gamma_{n-1}\end{array}\right)\right\}$$ where $\gamma_i$ stands
for a matrix in $\G_i=\text{GL}_i(R)$, $*$ for an arbitrary matrix
with entries in $R$, and $q^N*$ for a matrix of
the appropriate size with entries in $q^NR$, respectively. Then the set
$$\left\{(x,\bfy)\in P\times \G|\,v(x)=N, \nu(\mcR(\bfy^1 P^N))=\bfm,
\nu(\mcS(\bfy^1 P^N))=\bfn\right\} $$ may be written as a disjoint
union of the $|\G:\G_{\{1\},N}|$ sets
$$\left\{(x,\bfy)\in P\times \gamma\G_{\{1\},N}|\,v(x)=N,
\nu(\mcR(\bfy^1 P^N))=\bfm, \nu(\mcS(\bfy^1 P^N))=\bfn \right\}$$
where $\gamma$ runs through a complete set of coset representatives of
$\Gamma/\G_{\{1\},N}$. The measure of each of these sets is either
zero or equals $(1-q^{-1})q^{-N}\mu(\G_{\{1\},N})$. The latter happens
$\mcN'_{N,\bfm,\bfn}$ times, where $$\mcN'_{N,\bfm,\bfn}=\left|\left\{
\bfy^1\in\mathbb{P}^{n-1}(R/P^N)
|\;\ol{\bfy^1}\not=0,\nu(\mcR(\bfy^1))=\bfm,
\nu(\mcS(\bfy^1))=\bfn\right\}\right|.$$ Clearly
$\mcN_{N,\bfm,\bfn}=(1-q^{-1})q^N\mcN'_{N,\bfm,\bfn}$. Using the
identity $$\mu(\G)/\mu(\G_{\{1\},N})=\binom{n}{1}_{q^{-1}}q^{N(n-1)}$$
we obtain
$$ \mu_{N,\bfm,\bfn}=\mcN'_{N,\bfm,\bfn}(1-q^{-1})q^{-N}\mu(\G_{\{1\},N})
=\mcN_{N,\bfm,\bfn}\frac{\mu(\G)}{\binom{n}{1}_{q^{-1}}q^{N(n+1)}}$$
as claimed.
\end{proof}

Lemma~\ref{lemma relation eldiv} yields
\begin{align*}
P(\bfr,\bfs)-1
&=P(\bfr,\bfs)-\wt{Z_\varnothing}(-\bfr,-\bfs,\sum_\kappa r_\kappa +
\sum_\lambda s_\lambda -n-1)&&\eqref{Z empty rep growth}\\
&=\sum_{\substack{N\in\N \\\bfm\in\N_0^{k},
\bfn\in\N_0^{l}}}\mcN_{N,\bfm,\bfn}q^{-\sum_{\kappa\in[k]}(N-m_\kappa)r_\kappa-\sum_{\lambda\in[l]}(N-n_\lambda)s_\lambda}&&\eqref{P
rep growth} \\
&=\frac{\binom{n}{1}_{q^{-1}}}{\mu(\G)}\cdot\\&\hspace{0.2in}\sum_{\substack{N\in\N\\\bfm\in\N_0^f,\bfn\in\N_0^h}}\mu_{N,\bfm,\bfn}q^{N(n+1-\sum_{\kappa}r_\kappa-\sum_{\lambda}s_\lambda)+\sum_{\kappa}m_\kappa
r_\kappa + \sum_{\lambda}n_\lambda s_\lambda}&&\eqref{relation N mu
rep growth}\\
&=\frac{\binom{n}{1}_{q^{-1}}}{\mu(\G)}Z_{\{1\}}(-\bfr,-\bfs,\sum_\kappa
r_\kappa + \sum_\lambda s_\lambda -n-1)&&\eqref{Z 1 rep growth}\\
&=(1-q^{-n})\wt{Z_{\{1\}}}(-\bfr,-\bfs,\sum_\kappa r_\kappa +
\sum_\lambda s_\lambda -n-1).&&\eqref{definition Z I tilde}
\end{align*}

From Corollary~\ref{corollary 1 to inversion property} we
deduce\m{proposition elementary divisors}

\begin{proposition}\label{proposition elementary divisors}
For all but finitely many completions $K$ of $F$, the following
functional equation holds:
$$P_{\mcR,\mcS,K}(\bfr,\bfs)|_{q\rightarrow q^{-1}}=q^n
P_{\mcR,\mcS,K}(\bfr,\bfs).$$
\end{proposition}


\section{Applications to zeta functions of groups and rings}\label{section applications}
\subsection{Zeta functions of rings}\label{subsection rings}

In this section we prove Theorem~\ref{theorem rings}. Let $L$ be a
ring of torsion-free rank $n$. In fact, without loss of generality we
may assume that $L$ is additively isomorphic to $\Z^n$. Let $p$ be a
prime. Multiplication in $L$ is a bi-additive mapping
$\beta:\Z^n\times\Z^n\rightarrow\Z^n$, which extends to
$\beta_p:\Zp^n\times\Zp^n\rightarrow\Zp^n$, inducing a $\Zp$-algebra
structure on $L_p:=L\otimes\Zp$. We shall give a formula for the local
zeta functions
$$\zeta_{L,p}(s)=\sum_{H\leq L_p} |L_p:H|^{-s},$$ where $H$ runs
over the subalgebras of finite index in~$L_p$, valid for almost all
primes~$p$, in terms of the $p$-adic integrals studied in
Section~\ref{section igusa}. More precisely, we shall show that
$\zeta_{L,p}(s)$ is expressible in terms of functions $\wt{Z}(\bfs)$,
defined as in~\eqref{definition Z tilde}, to which
Corollary~\ref{corollary 2 to inversion property} is applicable.

Write $L=\Z l_1\oplus\dots\oplus\Z l_n$. We
consider the $n\times n$-matrix of $\Z$-linear forms
$$\mcR(\bfy)=(L_{ij}(\bfy))\in\text{Mat}_n(\Z[\bfy]),$$ where
$L_{ij}(\bfy):=\sum_{k\in[n]}\lambda_{ij}^ky_k$, encoding the
structure constants $\lambda_{ij}^k$ of $L$ with respect to the chosen
basis, that is $l_il_j=\sum_{k\in[n]}\lambda_{ij}^kl_k$. Let $C_i$
denote the matrix of the linear map given by right-multiplication with
the generator~$l_i$.

A full sublattice $\L$ in $(L_p,+)$ corresponds to a coset $\G M$,
where $\G=\G_n=\GL_n(\Zp)$ and the rows of the matrix
$M=(m_{ij})\in\text{GL}_n(\Qp)\cap\text{Mat}_n(\Zp)$, the set of
integral $n\times n$-matrices with non-zero determinant, encode the
coordinates of generators for $\L$ with respect to the chosen
basis. Denote by $M_i$ the $i$-th row of~$M$. It is not hard to check
(cf.~the proof of \cite[Theorem~5.5]{duSG/00}) that $\Lambda$ is a
$\Zp$-sub\emph{algebra} of $L_p$ if and only if\m{subalgebra
condition}
\begin{equation}\label{subalgebra condition}
\forall i,j\in[n]:\; M_i\sum_{r\in[n]} C_r m_{jr} \in\langle
M_k|\;k\in[n]\rangle_{\Zp}.
\end{equation}
Rather than trying to analyse the restrictions
condition~\eqref{subalgebra condition} imposes on the entries of
suitable upper-triangular representatives of the coset $\G M$ as
in~\cite{duSG/00}, we base our analysis on the following two basic
observations.

The first point is that every homothety class of lattices $\L$ in
$(L_p,+)$ contains a largest subalgebra $\L_0$, and the subalgebras in
this class are exactly the multiples $p^m\L_0$, $m\in\N_0$. We thus
have
\begin{equation}\label{reduction}
\zeta_{L,p}(s)=(1-p^{-ns})^{-1}\sum_{[\L]}|L_p:\L_0|^{-s},
\end{equation}
where $\L_0$ denotes the largest subalgebra in the homothety
class~$[\L]=[\L_0]$.

The second observation is that it is easy to check
condition~\eqref{subalgebra condition} if $M$ happens to be a diagonal
matrix. With respect to the given basis, $\G M$ may not admit a
diagonal representative. By the elementary divisor theorem, however,
it does contain a representative of the form $M=D\alpha^{-1}$, where $\alpha\in\G$,
  $$D=D(I,\bfr_0)=p^{r_0}\,\text{diag}(\underbrace{\underbrace{p^{\sum_{\iota\in
        I}r_\iota},\dots,p^{\sum_{\iota\in
        I}r_\iota}}_{i_1},\dots,p^{r_{i_l}},\dots,p^{r_{i_l}}}_{i_l},1,\dots,1)$$
        for a set $I=\{i_1,\dots,i_l\}_<\subseteq[n-1]$ and a vector
        $(r_0,r_{i_1},\dots,r_{i_l})=:\bfr_0\in\N_0\times\N^l$ (both
        depending only on $\G M$). We say that $\Lambda$ has type
        $(I,\bfr_0)$. We call $\L$ \emph{maximal} (in its homothety
        class) if $r_0=0$. We say that the homothety class $[\L]$ has
        type $(I,\bfr)$, $\bfr=(r_{i_1},\dots,r_{i_l})\in\N^l$, --
        written $\nu([\L])=(I,\bfr)$ -- if its maximal element has
        type $(I,(0,r_{i_1},\dots,r_{i_l}))$. By slight abuse of
        notation we may also say that a lattice $\L$ has type $I$ if
        $\nu([\L])=(I,\bfr)$ for some integral vector $\bfr$, and that
        a homothety class $[\L]$ has type $I$ if any of its elements
        does. In this case we write $\nu([\L])=I$. We shall denote by
        $\alpha^j$ the $j$-th column of the matrix $\alpha$ and by
        $D_{ii}$ the $i$-th diagonal entry of~$D$. Note that $\alpha$
        is only unique up to right-multiplication by an element of
$$\G_{I,\bfr}:=\text{Stab}_\G(\G D)=\left\{\left(\begin{tabular}{c|c|c|c|c}
$\gamma_{i_1}$&$*$&$\hdots$&$*$&$*$\\\hline
$ p^{r_{i_1}}*$ &$\gamma_{i_2-i_1}$&$\ddots$&$\vdots$&\\\hline
$p^{r_{i_1}+r_{i_2}}*$&$p^{r_{i_2}}*$&$\ddots$&$*$&$\vdots$\\\hline
$\vdots$&$\vdots$&&$\gamma_{i_l-i_{l-1}}$&$*$\\\hline
$p^{r_{i_1}+\dots+r_{i_l}}*$&$p^{r_{i_2}+\dots+r_{i_l}}*$&$\hdots$&$p^{r_{i_l}}*$&$\gamma_{n-i_l}$
\end{tabular}\right)\right\}$$
where $\gamma_\iota$ stands for a matrix in $\Gamma_\iota$, $*$ for an
arbitrary matrix with entries in $\Zp$, and $p^r*$ for a matrix with
entries in $p^r\Zp$ of the appropriate sizes, respectively. Thus there
is a $1-1$-correspondence between lattice classes $[\L]$ of type
$(I,(r_{i_1},\dots,r_{i_l}))$ and cosets
$\alpha\G_{I,\bfr}$. Furthermore\m{equation index}
\begin{equation}\label{equation index}
\left|\left\{[\L]|\;\nu([\L])=(I,\bfr)\right\}\right|=|\G:\G_{I,\bfr}|=\mu(\G)/\mu(\G_{I,\bfr})=\binom{n}{I}_{p^{-1}}p^{\sum_{\iota\in
    I}r_\iota \iota(n-\iota)},
\end{equation} 
where, as usual, $\mu$ denotes the Haar measure on $\G$ normalised so
that $\mu(\G)=(1-p^{-1})\dots(1-p^{-n})$. As we have noted above it
coincides with the additive Haar measure on
$\text{Mat}_n(\Zp)\cong\Zp^{n^2}$, normalised so that
$\mu(\text{Mat}_n(\Zp))=1$.

It is now straightforward to check that \eqref{subalgebra condition}
is equivalent to\m{subalgebra conditionII}
\begin{equation}\label{subalgebra conditionII}
\forall i\in[n]:\; D \mcR_{(i)}^<(\alpha) D \equiv 0 \mod D_{ii},
\end{equation}
where $\mcR_{(i)}^<(\alpha):=\alpha^{-1}\mcR(\alpha^i)(\alpha^{-1})^t.$ It
is easy to verify that condition \eqref{subalgebra conditionII} is
equivalent to
$$\forall
  i,r,s\in[n]:\;(\mcR_{(i)}^<(\alpha))_{rs}\,p^{r_0+\sum_{s\leq\iota\in
  I}r_\iota + \sum_{r\leq\iota\in I}r_\iota + \sum_{i>\iota\in
  I}r_\iota}\equiv 0 \mod p^{\sum_{\iota\in
  I}r_\iota}\label{subalgebra contitionIII}$$ which may in turn be
  reformulated as\m{subalgebra conditionIV}
\begin{equation}\label{subalgebra conditionIV}
r_0\geq\sum_{\iota\in I}r_\iota-\underbrace{\min\left\{\sum_{\iota\in
I}r_\iota,\sum_{s\leq\iota\in I}r_\iota + \sum_{r\leq\iota\in I}r_\iota +
\sum_{i>\iota\in
I}r_\iota+v_{irs}(\alpha)|\,(i,r,s)\in[n]^3\right\}}_{=:m([\L])},
\end{equation}
where
$v_{irs}(\alpha):=\min\left\{v_p\left((\mcR_{(\iota)}^<(\alpha))_{\rho
\sigma}\right)|\iota\leq i,\rho\geq r,\sigma\geq s\right\}$.

\begin{remark}\label{remark rings} 
Whilst it might seem more natural to replace the inequalities in this
definition of $v_{irs}(\alpha)$ by equalities, the present formulation
is preferable as it allows us to translate the counting problem into
the language developed in Section~\ref{section igusa}.
\end{remark}
Note that the right hand side of \eqref{subalgebra conditionIV}
depends only on the homothety class of $\L$ and is neatly separated in
terms that depend on the type $(I,\bfr)$ of $[\L]$ and terms that
depend solely on $\alpha$. The subalgebra $\L_0$ is characterised by
equality in~\eqref{subalgebra conditionIV}. As we shall see, this
formulation of the `subalgebra condition'~\eqref{subalgebra condition}
therefore enables us to express the $p$-th local zeta function of $L$
in terms of $p$-adic integrals on $(p\Zp)^{|I|}\times \G$,
$I\subseteq[n-1]$.

By \eqref{reduction} it suffices to compute
\begin{equation}\label{A}
A^<(s) := \sum_{[\L]}|L_p:\L_0|^{-s} = \sum_{I\subseteq[n-1]}
\underbrace{\sum_{\nu([\L])=I}|L_p:\L_0|^{-s}}_{=:A^<_I(s)}.
\end{equation}
For a fixed $I=\{i_1,\dots,i_l\}_<\subseteq[n-1]$, we set
\begin{equation*}
\mcN^<_{I,\bfr,m}:=|\{[\L]|\;\nu([\Lambda])=(I,\bfr),
 m([\L])=m\}|
\end{equation*} 
with $m([\Lambda])$ defined as in \eqref{subalgebra conditionIV}. Thus
\begin{align}
A^<_I(s)=&\sum_{\bfr=(r_{i_1},\dots,r_{i_l})\in\N^l}p^{-s\sum_{\iota\in
    I}\iota
    r_\iota}\sum_{m\in\N_0}\mcN^<_{I,\bfr,m}p^{-sn\left(\sum_{\iota\in
    I}r_\iota-m\right)}\nonumber\\
    =&\sum_{\bfr\in\N^l}p^{-s\sum_{\iota\in
    I}r_\iota(\iota+n)}\sum_{m\in\N_0}\mcN^<_{I,\bfr,m}p^{snm}.\label{A_I}
\end{align}

As in the proof of Proposition~\ref{proposition elementary divisors}
we shall design a $p$-adic integral to describe the generating
functions~$A^<_I(s)$. Consider the $p$-adic integral
$Z_I(\bfs)=Z_{W,\Qp,I}(s_1,\dots,s_n)$ defined as in
Section~\ref{section igusa} with $k=n$, $m=n^2$, $W=\G$ and set
\begin{align}
\bff_{irs}(\bfy)=&\{(\mcR_{(\iota)}^<(\bfy))_{\rho\sigma}|\;\iota\leq
  i,\rho\geq r,\sigma\geq s\}, (i,r,s)\in[n]^3,\nonumber\\
  \bfg_{n,I}(\bfx,\bfy)=&\left\{\prod_{\iota\in
  I}x_\iota\right\}\cup\bigcup_{(i,r,s)\in[n]^3}\left(\prod_{\iota\in
  I}x_\iota^{\delta_{\iota\geq r}+\delta_{\iota\geq
  s}+\delta_{\iota<i}}\right)\bff_{irs}(\bfy),\nonumber\\
  \bfg_{\kappa,I}(\bfx,\bfy) =& \left\{\prod_{\iota\in
  I}x_\iota^{\delta_{\iota\kappa}}\right\},\;\kappa\in[n-1].\label{f_irs}
\end{align}
The ideals $(\bff_{irs}(\bfy))$ can easily be seen to be
$B(\Qp)$-invariant. Without loss of generality we may assume that none
of them equals the zero ideal (otherwise we just omit the respective
ideal), and, by omitting at most finitely many primes, we may assume
that this also holds$\mod p$. (Note that we may well have to omit all
the $(\bff_{irs}(\bfy))$. This happens when the ring structure is
trivial.) Thus Theorems~\ref{theorem denef projective}
and~\ref{inversion property Z tilde} and their Corollaries apply to a
principalisation for the ideal $\mcI=\prod_{irs}(\bff_{irs}(\bfy))$
with good reduction$\mod p$. Note that
\begin{equation}\label{Z_I}
Z^<_I(\bfs):=Z_I(\bfs)=Z_I((s_\iota)_{\iota\in I},s_n)=\sum_{\bfr\in\N^l}p^{-\sum_{\iota\in I}s_\iota r_\iota}\sum_{m\in\N_0}\mu^<_{I,\bfr,m}p^{-s_nm},
\end{equation}
where 
\begin{equation*}
\mu^<_{I,\bfr,m}=\mu\left\{(\bfx,\bfy)\in (p\Zp)^l\times
\G|\,v_p(x_\iota)=r_\iota, m(\bfx,\bfy)=m\right\},
\end{equation*}
with
\begin{multline*}
m(\bfx,\bfy):=\\
 \min\left\{\sum_{\iota\in I}r_\iota,\sum_{s\leq\iota\in I}v_p(x_\iota) +
\sum_{r\leq\iota\in I}v_p(x_\iota) + \sum_{i>\iota\in
I}v_p(x_\iota)+v_{irs}(\bfy)|\;(i,r,s)\in[n]^3\right\}.
\end{multline*}

Again we need a lemma to relate the numbers $\mu^<_{I,\bfr,m}$ to the
data $\mcN^<_{I,\bfr,m}$ we are trying to understand.\m{lemma relation
rings}

\begin{lemma}\label{lemma relation rings}
\m{relation N mu}
\begin{equation}\label{relation N mu}
\mcN^<_{I,\bfr,m}=\frac{\mu^<_{I,\bfr,m}}{(1-p^{-1})^lp^{-\sum_{\iota\in
      I}r_\iota}\mu(\G_{I,\bfr})}=\frac{\binom{n}{I}_{p^{-1}}}{(1-p^{-1})^l\mu(\G)}\mu^<_{I,\bfr,m}\,p^{\sum_{\iota\in
      I}r_\iota(\iota(n-\iota)+1)}.
\end{equation}
\end{lemma}

\begin{proof} The set $$\left\{(\bfx,\bfy)\in (p\Zp)^l\times
\G|\,v_p(x_\iota)=r_\iota, m(\bfx,\bfy)=m\right\}$$ may be written as
a disjoint union of the $|\G:\G_{I,\bfr}|$ sets
$$\left\{(\bfx,\bfy)\in (p\Zp)^l\times
\gamma\G_{I,\bfr}|\,v_p(x_\iota)=r_\iota, m(\bfx,\bfy)=m\right\}$$
where $\gamma$ runs through a complete set of coset representatives
for~$\G/\G_{I,\bfr}$. The measure of each of these sets is either zero
or equal to $(1-p^{-1})^lp^{-\sum_{\iota\in
I}r_\iota}\mu(\G_{I,\bfr})$. The latter happens if and only if
$\gamma\G_{I,\bfr}$ corresponds to a lattice of type $(I,\bfr)$ such
that $m([\L])=m$. This proves the first equality. The second equality
is an immediate consequence of~\eqref{equation index}.
\end{proof} 

Lemma~\ref{lemma relation rings} allows us to express the generating
functions~\eqref{A_I} in terms of the $p$-adic
integrals~\eqref{Z_I}. Indeed,

\begin{align*}
A^<_I(s)&=\sum_{\bfr\in\N^l}p^{-s\sum_{\iota\in
      I}r_\iota(\iota+n)}\sum_{m\in\N_0}\mcN^<_{I,\bfr,m}p^{snm}&&\eqref{A_I}\\
      &=\frac{\binom{n}{I}_{p^{-1}}}{(1-p^{-1})^l\mu(\G)}\sum_{\bfr}p^{-\sum_{\iota\in
      I}r_\iota(s(\iota+n)-\iota(n-\iota)-1)}\sum_m\mu^<_{I,\bfr,m}p^{snm}&&\eqref{relation
      N mu}\\
      &=\frac{\binom{n}{I}_{p^{-1}}}{(1-p^{-1})^l\mu(\G)}Z^<_I((s(\iota+n)-\iota(n-\iota)-1)_{\iota\in
      I},-sn)&&\eqref{Z_I}\\
      &=\binom{n}{I}_{p^{-1}}\wt{Z^<_I}((s(\iota+n)-\iota(n-\iota)-1)_{\iota\in
      I},-sn).&&\eqref{definition Z I tilde}
\end{align*}
Theorem~\ref{theorem rings} follows now from equations~\eqref{reduction}
and \eqref{A}, and Corollary~\ref{corollary 2 to inversion property}.


\subsection{Conjugacy zeta functions of nilpotent groups}\label{subsection conjugacy}
In this section we prove Theorem~\ref{theorem conjugacy}. As mentioned in the introduction, it suffices to compute\m{definition conjugacy zeta function}
\begin{equation}\label{definition conjugacy zeta function}
\zetacc_{L,p}(s)=\sum_{H\leq
L_p}|L_p:H|^{-s}|L_p:\mathcal{N}_{L_p}(H)|^{-1}
\end{equation} 
for almost all primes $p$, where $L=L(G)$, $L_p=L\otimes \Zp$ and
$\mathcal{N}_{L_p}(H)=\{l\in L_p|\;[H,l]\leq H\}$ is the normaliser of $H$ in $L_p$.

\begin{remark} 
It was pointed out, e.g. in~\cite[Section 3.8]{Woodward/05}, that,
whilst~\eqref{definition conjugacy zeta function} might be used to
define local `conjugacy' zeta functions for arbitrary rings, they
might not encode the solutions to any actual counting problem unless
$L_p=L\otimes\Zp$ for a nilpotent Lie ring~$L$. It is for that reason
that we formulate Theorem~\ref{theorem conjugacy} in terms of
$\T$-groups $G$, bearing in mind that its proof is immediately reduced
to a computation in the Lie ring $L(G)$ that draws upon neither the
nilpotency nor the Lie property of $L(G)$.
\end{remark}
Along the lines of the `first observation' in the proof of
Theorem~\ref{theorem rings} we use the fact that, for a subring $H$ of
$L_p$, the normaliser $\mathcal{N}_{L_p}(H)$ is an invariant of the
homothety class~$[H]$ of $H$. We may thus write
$\mathcal{N}_{L_p}([H])$ for $\mathcal{N}_{L_p}(H)$. Resuming the
notation of the proof of Theorem~\ref{theorem rings} we have
\m{conjugacy chain}
\begin{equation}\label{conjugacy chain}
L_p\geq\mathcal{N}_{L_p}([H])\geq\L_0\geq p\L_0\geq p^2\L_0\geq\dots
\end{equation}
($H\in\{p^m\L_0|\;m\in\N_0\}$). In addition to recording the index
$|L_p:\L_0|$ as we run over homothety classes $[\L]$ of lattices in
$\Qp^n$, we now have to control the index
$|L_p:\mathcal{N}_{L_p}([\L])|$. We shall see that this index, too,
might be expressed in terms of congruences involving the type
$(I,\bfr)$ and coset $\alpha\G_{I,\bfr}$ determining $[\L]=[\L_0]$,
similar to the congruences~\eqref{subalgebra conditionII}. To compute
their index requires us to control the elementary divisors of certain
matrices. A similar problem was considered in Section~\ref{subsection
elementary divisors}; the problem we will have to solve in
Section~\ref{subsection normal} is also of this kind (though slightly
simpler). We resume the notation from Section~\ref{subsection
rings}. From \eqref{conjugacy chain} we deduce
\begin{equation*}
\zetacc_{L,p}(s)=(1-p^{-ns})^{-1}\sum_{[\L]}|L_p:\L_0|^{-s}|L_p:\mathcal{N}_{L_p}([\L])|^{-1}.
\end{equation*}
Using a fixed basis to identify $L_p$ with $\Zp^n$ we may express the
condition $\bfx\in\mathcal{N}_{L_p}([H])$ in a similar fashion
to~\eqref{subalgebra condition}. If $M$ is any matrix whose rows $M_k$
encode the coordinates of generators for any element of $[\L]$, a
vector $\bfx=(x_1,\dots,x_n)\in\Zp^n$ is in $\mathcal{N}_{L_p}([\L])$
if and only if \m{conjugacy condition}
\begin{equation}\label{conjugacy condition}
\forall i\in[n]: \;\bfx\sum_{r\in[n]}C_r
m_{ir}\in\langle M_k|\;k\in[n]\rangle_{\Zp}.
\end{equation}
Using the correspondence between lattice classes and pairs
$((I,\bfr),\alpha\G_{I,\bfr})$, condition \eqref{conjugacy condition}
may be reformulated as \m{conjugacy conditionII}
\begin{equation}\label{conjugacy conditionII}
\forall i\in[n]: \;\bfx \mcR_{(i)}^{\cc}(\alpha)D\equiv 0 \mod D_{ii},
\end{equation}
where $\mcR_{(i)}^\cc(\alpha):=\mcR(\alpha^i)(\alpha^{-1})^t$ and
$D=D(I,\bfr)$. We note that the scalar $p^{r_0}$ in $D$ cancels
in~\eqref{conjugacy conditionII}; we may thus
assume~$r_0=0$. Condition~\eqref{conjugacy conditionII} is then
equivalent to\m{conjugacy conditionIII}
\begin{equation}\label{conjugacy conditionIII}
\forall i\in[n]:\;\bfx\mcR_{(i)}^\cc(\alpha) D
p^{\sum_{\iota<i}r_\iota}\equiv 0 \mod p^{\sum_{\iota\in I}r_\iota}.
\end{equation}
Setting
\begin{align*}
\mcR^\cc(\alpha)&=(\mcR_{(1)}^\cc(\alpha)|\dots|\mcR_{(n)}^\cc(\alpha)),\\
\bfD^\cc(I,\bfr)&=\text{diag}(\underbrace{D,\dots,D}_{i_1},\underbrace{p^{r_{i_1}}D,\dots,p^{r_{i_1}}D}_{i_2-i_1},\dots,\underbrace{p^{\sum_{\iota\in
I}r_\iota}D,\dots,p^{\sum_{\iota\in I}r_\iota}D}_{n-i_l}),
\end{align*}
(where `$\text{diag}$' refers to the diagonal $n^2\times n^2$-matrix
built from $n$ scalar multiples of the diagonal $n\times n$-matrix~$D$),
\eqref{conjugacy conditionIII} may in turn be reformulated
as\m{conjugacy conditionIV}

\begin{equation}\label{conjugacy
conditionIV}
\bfx\mcR^\cc(\alpha)\bfD^\cc(I,\bfr)\equiv 0\mod p^{\sum_{\iota\in I}r_\iota}
\end{equation}
To keep track of the index of $\mcN_{L_p}([\L])$, the full sublattice
of $L_p\cong\Zp^n$ of solutions to~\eqref{conjugacy conditionIV}, we
introduce an invariant $\nu^\cc([\Lambda])\in\N_0^n$ as follows. We
say that $\nu^\cc([\L])=\bfm=(m_1,\dots,m_n)\in\N_0^n$ if
\begin{itemize}
\item the matrix $\mcR^\cc(\alpha)\bfD^\cc(I,\bfr)$ has elementary
divisor type $$\bfwtm=(\wt{m}_1,\dots,\wt{m}_n)$$ (i.e. there are
matrices $\beta\in\G_n$, $\gamma\in\G_{n^2}$ such that
$\beta\mcR^\cc(\alpha)\bfD^\cc(I,\bfr)\gamma=\left(\text{diag}(p^{\wt{m}_1},\dots,p^{\wt{m}_n})|0\right)$
and $\wt{m}_i\in\N_0\cup\{\infty\}$, $\wt{m}_1\leq\dots\leq\wt{m}_n$)
and
\item  $\bfm=(m_1,\dots,m_n)$ is defined by $m_i=\min\left\{\sum_{\iota\in I}r_\iota,\wt{m}_i\right\}$ for each $i\in[n]$.
\end{itemize}
The index $|L_p:\mcN_{L_p}([\L])|$ equals
$$p^{\sum_{i\in[n]}\left(\sum_{\iota\in I}r_\iota-m_i\right)}.$$
As in Section~\ref{subsection rings} it is helpful to write
\begin{align*}
A^\cc(s):=&\sum_{[\L]}|L_p:\L_0|^{-s}|L_p:\mathcal{N}_{L_p}([\L])|^{-1}\\
=&\sum_{I\subseteq[n-1]}\underbrace{\sum_{\nu([\L])=I}|L_p:\L_0|^{-s}|L_p:\mathcal{N}_{L_p}([\L])|^{-1}}_{=:A^\cc_I(s)}.
\end{align*}
For a fixed $I=\{i_1,\dots,i_l\}_<\subseteq[n-1]$, we set
$$\mcN^\cc_{I,\bfr,m,\bfm}:=|\{[\L]|\;\nu([\L])=(I,\bfr),m([\L])=m,\nu^\cc([\L])=\bfm\}|$$
with $m([\L])$ defined as in~\eqref{subalgebra conditionIV}. We thus
have\m{A cc}
\begin{align}
A^\cc_I(s)\nonumber\\ =&\sum_{\bfr\in\N^l}p^{-s\sum_{\iota\in I}\iota
r_\iota}\sum_{\substack{m\in\N_0\\\bfm\in\N_0^n}}\mcN^\cc_{I,\bfr,m,\bfm}p^{-sn(\sum_{\iota\in
I}r_\iota-m)-\sum_{i\in[n]}(\sum_{\iota\in I}r_\iota-m_i)}\nonumber\\
=&\sum_{\bfr\in\N^l}p^{\sum_{\iota\in
I}r_\iota(-s(\iota+n)-n)}\sum_{\substack{m\in\N_0\\\bfm\in\N_0^n}}\mcN^\cc_{I,\bfr,m,\bfm}p^{snm+\sum_{i\in[n]}m_i}.\label{A
cc}
\end{align}
As in Section~\ref{subsection rings} we shall show that the generating
functions $A^\cc_I(s)$ may be expressed in terms of the $p$-adic
integrals to which the results of Section~\ref{section igusa} may be
applied. In order to control the invariant $\nu^{\cc}([\Lambda])$ we
need to parametrise the minors of the matrix $\mcR^\cc(\bfy)$ in a
suitable way. This motivates the following combinatorial definitions.

Let $\matnzeroone$ denote the set of $n\times n$-matrices with entries
in $\{0,1\}$, and
$\mcS^\cc_{j,n}=\{S\in\matnzeroone|\sum_{r,s}S_{rs}=j\}$. We introduce
a partial order on $\mcS^\cc_{j,n}$ by saying that, given
$S,T\in\mcS^\cc_{j,n}$, $T\preceq S$ if,
 
\begin{tabular}{rl}for all $r\in[n]$:&
$\sum_{\rho\leq r}\sum_{\sigma\in[n]}T_{\rho\sigma}\geq\sum_{\rho\leq
r}\sum_{\sigma\in[n]}S_{\rho\sigma}$ and, \\for all $s\in[n]$:&
$\sum_{\sigma>s}\sum_{\rho\in[n]}T_{\rho\sigma}\geq
\sum_{\sigma>s}\sum_{\rho\in[n]}S_{\rho\sigma}$.
\end{tabular}

\noindent Pictorially speaking, this amounts to saying that the matrix
$T$ may be obtained from the matrix $S$ by moving some of the non-zero
entries towards `north-east'. Given a matrix $S\in\mcS^\cc_{j,n}$ and a
$n\times n^2$-matrix $M=(M_1|\dots|M_n)$, $M_i\in\text{Mat}_n(\Zp)$, a
$j\times j$-submatrix of $M$ of column-type $S$ is a submatrix
obtained by choosing $j$ rows of $M$ and the $s$-th column of $M_r$ if
and only if $S_{rs}=1$. For $S\in\mcS^\cc_{j,n}$, we set \m{definition
fccjs}
\begin{equation}\label{definition fccjs}
\bff^\cc_{j,S}(\bfy)=\{\det(N)|\;N\text{ a $j\times j$-submatrix of
$\mcR^\cc(\bfy)$ of column-type $T\preceq S$}\}
\end{equation}
and define the monomial
\begin{equation*}
M^\cc_{S,I}(\bfx)=\prod_{\iota\in
I}x_\iota^{\sum_{\rho>\iota}\sum_{\sigma=1}^nS_{\rho\sigma}+\sum_{\sigma\leq\iota}\sum_{\rho=1}^nS_{\rho\sigma}}.
\end{equation*}
We set
$$\bff^\cc_{j,I}(\bfx,\bfy)=\bigcup_{S\in\mcS^\cc_{j,n}}M^\cc_{S,I}(\bfx)\bff^\cc_{j,S}(\bfy).$$
Now we define, using the sets of polynomials
$\bfg_{\kappa,I}(\bfx,\bfy)$ introduced in~\eqref{f_irs},\m{definition
Z I cc}
\begin{align}\label{definition Z I cc}
\lefteqn{Z^\cc_I(\bfs,\wt{\bfs},\wt{\wt{\bfs}})=\int_{(p\Zp)^l\times
W}\prod_{\kappa\in[n]}\|\bfg_{\kappa,I}(\bfx,\bfy)\|^{s_\kappa}\cdot}\nonumber\\&&\quad\prod_{j\in[n]}\left(\|\bff^\cc_{j,I}(\bfx,\bfy)\cup\left(\prod_{\iota\in
I}x_\iota\right)\bff^\cc_{j-1,I}(\bfx,\bfy)\|^{\wt{s}_j}\|\bff^\cc_{j-1,I}(\bfx,\bfy)\|^{\wt{\wt{s}}_j}\right)|\tud\bfx_I||\tud\bfy|.
\end{align}

\begin{remark} 
Note that $M^\cc_{T,I}(\bfx)|M^\cc_{S,I}(\bfx)$ if $T\preceq S$. We
therefore could have kept definition \eqref{definition fccjs} simpler
by replacing `$T\preceq S$' by `$S$', without changing the
integral~$Z^\cc_I$. The extra complication ensures that the results
from Section~\ref{section igusa} are applicable.
\end{remark}

We leave it to the reader to verify that, for each $j\in [n]$ and
$S\in \mcS^\cc_{j,n}$, the ideal $(\bff^\cc_{j,S}(\bfy))$ is
$B(\Qp)$-invariant. Therefore the Theorems~\ref{theorem denef
projective} and~\ref{inversion property Z tilde} and their corollaries
are applicable to the
integral~$Z^\cc_I(\bfs,\wt{\bfs},\wt{\wt{\bfs}})$. Note that
$Z^\cc_I(\bfs,\wt{\bfs},\wt{\wt{\bfs}})=Z^\cc_I((s_\iota)_{\iota\in
I},s_n,\wt{\bfs},\wt{\wt{\bfs}})$. We set
$$\mu^\cc_{I,\bfr,m,\bfm}=\mu\{(\bfx,\bfy)\in(p\Zp)^l\times
W|\nu_p(\bfx_\iota)=r_\iota,
m(\bfx,\bfy)=m,\nu^\cc(\bfx,\bfy)=\bfm\}.$$ Here
$\nu^\cc(\bfx,\bfy)=\bfm=(m_1,\dots,m_n)$ if, for each $i\in[n]$,
$$m_i=\min\left\{\sum_{\iota\in I}v_p(x_\iota),\wt{m}_i\right\},$$
where $\wt{\bfm}=(\wt{m}_1,\dots,\wt{m}_n)$ is the elementary divisor
type of the matrix $$\mcR^\cc(\bfy)\bfD^\cc(I,(v_p(x_\iota))_{\iota\in
I}).$$ Then
\begin{equation*}
\begin{split}
Z^\cc_I((s_\iota)_{\iota\in
I},s_n,\wt{\bfs})&:=Z^\cc_I((s_\iota)_{\iota\in
I},s_n,\wt{\bfs},-{\wt{\bfs}})\\&=\sum_{\bfr\in\N^l}p^{-\sum_{\iota\in
I}s_\iota
r_\iota}\sum_{\substack{m\in\N_0\\\bfm\in\N^n}}\mu^\cc_{I,\bfr,m,\bfm}p^{-s_nm-\sum_{i\in[n]}\wt{s}_im_i}.
\end{split}
\end{equation*}
As in~\eqref{relation N mu} we want to relate the numbers
$\mu^\cc_{I,\bfr,m,\bfm}$ with the integers
$\mcN^\cc_{I,\bfr,m,\bfm}$.

\begin{lemma}\label{lemma relation cc}
\m{relation N mu cc}
\begin{equation}\label{relation N mu cc}
  \mcN^\cc_{I,\bfr,m,\bfm}=\frac{\binom{n}{I}_{p^{-1}}}{(1-p^{-1})^l\mu(\G)}\mu^\cc_{I,\bfr,m,\bfm}p^{\sum_{\iota\in
  I}r_\iota(\iota(n-\iota)+1)}.
\end{equation}
\end{lemma}

\begin{proof} Analogous to the proof of Lemma~\ref{lemma relation rings}.
\end{proof}

Thus
\begin{align*}
A^\cc_I(s)\\ =&\sum_{\bfr\in\N^l}p^{\sum_{\iota\in
I}r_\iota(-s(\iota+n)-n)}\sum_{\substack{m\in\N_0,\\\bfm\in\N_0^n}}\mcN^\cc_{I,\bfr,m,\bfm}p^{snm+\sum_{i\in[n]}m_i}&&\eqref{A
cc}\\ =&\frac{\binom{n}{I}_{p^{-1}}}{(1-p^{-1})^l\mu(\G)}\cdot\\
&\hspace{0.5in}\sum_{\bfr\in\N^l}p^{\sum_{\iota}r_\iota(-s(\iota+n)+\iota(n-\iota)+1-n)}\sum_{m,\bfm}\mu^\cc_{I,\bfr,m,\bfm}p^{snm+\sum_{i\in[n]}m_i}&&\eqref{relation
N mu cc}\\ =&\frac{\binom{n}{I}_{p^{-1}}}{(1-p^{-1})^l\mu(\G)}\cdot\\
&\hspace{0.5in}{Z^\cc_I}((s(\iota+n)-\iota(n-\iota)-1+n)_{\iota\in
I},-sn,-1,\dots,-1)&&\eqref{definition Z I cc}\\
=&\binom{n}{I}_{p^{-1}}\wt{Z^\cc_I}((s(\iota+n)-\iota(n-\iota)-1+n)_{\iota\in
I},-sn,-1,\dots,-1).&&\eqref{definition Z I tilde}
\end{align*}

\noindent Theorem~\ref{theorem conjugacy} now follows from
Corollary~\ref{corollary 2 to inversion property}.


\subsection{Normal zeta functions of class-$2$-nilpotent groups}\label{subsection normal}
In this section we prove Theorem~\ref{theorem normal}. Let $G$ be a
$\Ttwo$-group as in the statement of the theorem. Without loss of
generality we may assume that $G/Z(G)$ and $Z(G)$ are torsion-free
abelian groups of rank $d$ and~$d'$, respectively (so $n=d+d'$), and
that $G$ admits a presentation
$$G=\langle
g_1,\dots,g_{d},h_1,\dots,h_{d'}|\,[g_i,g_j]=\sum_{k\in[d]}\lambda_{ij}^kh_k,
\,\lambda_{ij}^k\in\Z,\text{ all other $[\,,]$ trivial}\rangle.$$
(Note that we used additive notation for expressions in the abelian
group $G'$.) Thus we obtain a matrix
$$\mcM(\bfy):=\left({L}_{ij}(\bfy)\right)\in
\text{Mat}_d(\Z[\bfy])$$ of linear forms
${L}_{ij}(\bfy):=\sum_{k\in[d]}\lambda_{ij}^ky_k$, encoding the
commutator structure of $G$.  By disregarding at most finitely many
further primes we may also assume that $p$ does not divide
$\mcM(\alpha)$ whenever $\alpha\in\Zp^{d'}\setminus p\Zp^{d'}$.

We begin our argument as in \cite[Section 3]{Voll/05}, albeit with
slightly different notation. Note, however, that we do not require the
assumption that $Z(G)=G'$. By \cite[Lemma 1]{Voll/05} and its
corollary, it suffices to compute a functional equation for the
generating function
$$A^\tl(s)=\sum_{I\subseteq[d'-1]}A^\tl_I(s),$$ where
$$A^\tl_I(s)=\sum_{\nu([\L])=I}|Z(L_p):\L|^{d-s}|L_p:X(\L)|^{-s}.$$
Here the sum ranges over homothety classes of \emph{maximal} lattices
$\L$ of type\footnote{Note that the definition of a lattice's type
given in~\cite{Voll/05} differs from the one in the current paper in
so far as $I$ is replaced by $d'-I=\{d'-i|\,i\in I\}$.}~$I$ in the
centre $Z(L_p)$ of the $\Zp$-algebra
$$L_p:=(G/Z(G)\oplus Z(G))\otimes\Zp,$$ and $X(\L)/\L=Z(L_p/\L)$. We
identify $Z(L_p)$ with $\Zp^{d'}$ using the basis
$\{h_1,\dots,h_{d'}\}$. The index $|L_p:X(\L)|$ is the index in
$L_p/Z(L_p)\cong\Zp^d$ of a system of linear congruences which we now
describe (cf.~\cite[\S 2.2]{Voll/04}). Let $[\L]$ be of type
$(I,\bfr)$, $I=\{i_1,\dots,i_l\}_<\subseteq[d'-1]$,
$\bfr=(r_{i_1},\dots,r_{i_l})\in\N^l$, corresponding to the coset
$\alpha\G_{I,\bfr}\in \G_{d'}/\G_{I,\bfr}$ as in the proof of
Theorem~\ref{theorem rings}. We set
\begin{align}
\mcM^\tl(\alpha)=&\left(\mcM(\alpha^1)|\dots|\mcM(\alpha^{d'})\right)\label{definition
M(alpha)},\\ \bfD^\tl(I,\bfr)=&\text{
diag}(\underbrace{1,\dots,1}_{di_1},\underbrace{p^{r_{i_1}},\dots,p^{r_{i_1}}}_{d(i_2-i_1)},\dots,\underbrace{p^{\sum_{\iota\in
I}r_\iota},\dots,p^{\sum_{\iota\in
I}r_\iota}}_{d(d'-i_{l})}).\nonumber
\end{align}
The system of linear congruences under consideration is
\begin{equation}
\bfx\mcM^\tl(\alpha)\bfD^\tl(I,\bfr)\equiv 0 \mod p^{\sum_{\iota\in
I}r_\iota}.\label{linear congruences}
\end{equation} The solutions to this system form a full lattice in~$L_p/Z(L_p)$. To keep track of its index we define the
invariant $\nu'([\L])\in\N_0^d$ as follows. We say that
$\nu'([\L])=\bfm=(m_1,\dots,m_d)$ if
\begin{itemize}
\item the matrix
$\mcM^\tl(\alpha)\bfD^\tl(I,\bfr)$ has elementary divisor type
$$\bfwtm=(\wt{m}_1,\dots,\wt{m}_d)$$ and 
\item $\bfm=(m_1,\dots,m_n)$ is defined by $m_i=\min\{\sum_{\iota\in
I}r_\iota,\wt{m}_i\}$ for each $i\in[n]$.
\end{itemize}
The index of the system~\eqref{linear congruences} equals
$$p^{\sum_{j\in[d]}\left(\sum_{\iota\in I}r_\iota - m_j\right)}.$$ By defining
$$\mcN^\tl_{I,\bfr,\bfm}:=|\{[\L]|\,\nu([\L])=(I,\bfr),\nu'([\L])=\bfm\}|$$
(note that $\mcN^\tl_{I,\bfr,\bfm}\not=0$ implies $m_1=0$, as
$p\not|\mcM^\tl(\alpha)$) we obtain a formula for $A^\tl_I(s)$ that is
analogous to \eqref{A_I}:
\begin{align}\label{A_I normal}
A^\tl_I(s)=&\sum_{\bfr\in\N^l}p^{(d-s)\sum_{\iota\in I}
r_\iota\iota}\sum_{\bfm\in\N_0^d}\mcN^\tl_{I,\bfr,\bfm}p^{-s\sum_{j\in[d]}(\sum_{\iota\in
I}r_\iota - m_j)}\nonumber\\
=&\sum_{\bfr}p^{\sum_{\iota\in I}r_{\iota}(-s(d+\iota)+\iota d)}\sum_{\bfm}\mcN^\tl_{I,\bfr,\bfm}p^{s\sum_{j\in[d]}m_j}.
\end{align}
As in Section~\ref{subsection rings} we shall show that the generating
function $A^\tl_I(s)$ may be expressed in terms of a $p$-adic integral
to which the results of Section~\ref{section igusa} can be applied. We
shall need some more notation.

Given the $d\times dd'$-matrix $M=(M_1|\dots|M_{d'})$,
$M_i\in\text{Mat}_d(\Zp)$, a $j\times j$-submatrix of $M$ of
column-type $S=(\sigma_1,\dots,\sigma_{d'})\in\N_0^{d'}$,
$\sum_{i\in[d']}\sigma_i=j$, is a submatrix obtained by choosing $j$
rows of $M$ and $\sigma_i$ columns in the `block' $M_i$ for each
$i\in[d']$. We denote by
$\mcS^\tl_{j,d'}=\{(\sigma_1,\dots,\sigma_{d'})|\sum_{i\in[d']}\sigma_i=j\}$
the set of possible such column-types. Given $S=(\sigma_i)$,
$T=(\tau_i)\in\mcS^\tl_{j,d'}$, we write $T\preceq S$ if, for all
$i\in[d']$, $\sum_{\iota\leq i}\tau_\iota\geq\sum_{\iota\leq
i}\sigma_\iota$. For $S\in\mcS^\tl_{j,d'}$, we set
\begin{equation}
\bff^\tl_{j,S}(\bfy)=\{\det(N)|\,N\text{ a $j\times j$-submatrix of
$\mcM(\bfy)$ of column-type }T\preceq S\}\label{definition bff S j}.
\end{equation}
 We define the monomial $M^\tl_{S,I}(\bfx)=\prod_{\iota\in
 I}x_\iota^{\sum_{\iota<\kappa\in[d']}\sigma_\kappa}$ and set
$$\bff^\tl_{j,I}(\bfx,\bfy)=\bigcup_{S\in\mcS^\tl_{j,d'}}M^\tl_{S,I}(\bfx)\bff^\tl_{j,S}(\bfy).$$
We are now ready to define
\begin{multline*}
Z_I^\tl(\bft,\bfs,\bfwts)=Z^\tl_I((t_\iota)_{\iota\in
I},s_2,\dots,s_d,\wt{s}_2,\dots,\wt{s}_d)=\int_{(p\Zp)^l\times
W} \prod_{\iota\in I}|x_\iota|^{t_\iota}\cdot\\
\quad\prod_{j\in[2,d]}\left(\|\bff^\tl_{j,I}(\bfx,\bfy)\cup\left(\prod_{\iota\in
I}x_\iota\right)\bff^\tl_{j-1,I}(\bfx,\bfy)\|^{s_j}\|\bff^\tl_{j-1,I}(\bfx,\bfy)\|^{\wt{s}_j}\right)|\tud\bfx_I||\tud\bfy|.
\end{multline*}

\begin{remark} 
Note that $M^\tl_{T,I}(\bfx)|M^\tl_{S,I}(\bfx)$ if $T\preceq S$. We
could have kept definition~\eqref{definition bff S j} simpler by
replacing `$T\preceq S$' by `$S$'. The extra complication ensures that
the results from Section~\ref{section igusa} are applicable. Note also
that we are not losing anything by omitting the factor for $j=1$, as
$\|\bff^\tl_{0,I}(\bfx,\bfy)\|=\|\bff^\tl_{1,I}(\bfx,\bfy)\|= 1$ for all
$\bfx\in(p\Zp)^l$, $\bfy\in W$.
\end{remark}

We leave it to the reader to verify that, for each $j\in[d]$ and
$S\in\mcS^\tl_{j,d'}$, the ideal $(\bff^\tl_{j,S}(\bfy))$ is
$B(\Qp)$-invariant. Therefore Theorems~\ref{theorem denef projective},
\ref{inversion property Z tilde} and their corollaries are
applicable. We set
$$\mu^\tl_{I,\bfr,\bfm}=\mu\{(\bfx,\bfy)\in (p\Zp)^l\times
W|\,v_p(x_\iota)=r_\iota,\nu'(\bfx,\bfy)=\bfm\},$$ where
$\nu'(\bfx,\bfy)=\bfm=(m_1,\dots,m_d)\in\N_0^d$ if
$\bfm=\min\{\sum_{\iota\in I}r_\iota,\bfwtm\}$ where $\bfwtm$ is the
elementary divisor type of the matrix
$\mcM^\tl(\bfy)\bfD^\tl(I,(v_p(x_\iota))_{\iota\in I})$. Note that
$\mu^\tl_{I,\bfr,\bfm}\not=0$ implies $m_1=0$. Then
\begin{equation}
Z^\tl_I((t_\iota)_{\iota\in I},\bfs):=Z^\tl_I((t_\iota)_{\iota\in
I},\bfs,-\bfs)=\sum_{\bfr\in\N^l}p^{-\sum_{\iota\in I}t_\iota
r_\iota}\sum_{\bfm\in\N_0^d}\mu^\tl_{I,\bfr,\bfm}p^{-\sum_{j\in[2,d]}s_j
m_j}.\label{Z_I normal}
\end{equation}\m{lemma relation normal}

\begin{lemma}\label{lemma relation normal}
\begin{equation}\label{relation N mu normal}
\mcN^\tl_{I,\bfr,\bfm}=\frac{\binom{d'}{I}_{p^{-1}}}{(1-p^{-1})^l\mu(\G)}\mu^\tl_{I,\bfr,\bfm}\,p^{\sum_{\iota\in
      I}r_\iota(\iota(d'-\iota)+1)}.
\end{equation}
\end{lemma}

\begin{proof} Analogous to the proof of Lemma~\ref{lemma relation rings}.
\end{proof}

\noindent Thus
\begin{align*}
A^\tl_I(s)\\=& \sum_{\bfr}p^{\sum_{\iota\in
I}r_{\iota}(-s(d+\iota)+\iota
d)}\sum_{\bfm}\mcN^\tl_{I,\bfr,\bfm}p^{s\sum_{j\in[d]}m_j}&&\eqref{A_I
normal}\\=&\frac{\binom{d'}{I}_{p^{-1}}}{(1-p^{-1})^l\mu(\G)}\cdot\\&\hspace{0.5in}\sum_{\bfr}p^{\sum_{\iota\in
I}r_\iota(-s(d+\iota)+\iota(d'+d-\iota)+1)}\sum_\bfm\mu^\tl_{I,\bfr,\bfm}(p^{-s})^{-\sum_{j\in[d]}m_j}&&\eqref{relation
N mu
normal}\\=&\frac{\binom{d'}{I}_{p^{-1}}}{(1-p^{-1})^l\mu(\G)}{Z_I}((s(d+\iota)-\iota(d+d'-\iota)-1)_{\iota\in
I},-s,\dots,-s)&&\eqref{Z_I
normal}\\=&\binom{d'}{I}_{p^{-1}}\wt{Z_I}((s(d+\iota)-\iota(d+d'-\iota)-1)_{\iota\in
I},-s,\dots,-s).&&\eqref{definition Z I tilde}
\end{align*}
Theorem~\ref{theorem normal} now follows from Corollary~\ref{corollary
2 to inversion property}.


\subsection{Representation zeta functions of $\T$-groups}\label{subsection howe}
In this section we recall some of Howe's work \cite{Howe-nilpotent/77}
on irreducible representations of $\T$-groups and co-adjoint orbits
and prove Theorem~\ref{theorem representations}.

If~$G$ is a group and $H\leq G$ is a subgroup we say that $H$ is
saturated if $g^n\in H$ implies $g\in H$ for all $g\in G$. Denote
by~$H_s$ the smallest saturated subgroup of~$G$ containing~$H$ (this
is called the {\sl isolator} of~$H$ in~\cite[Chapter~8,
Section~A]{Segal/83}). Clearly, if $H\triangleleft G$, $H=H_s$ if and
only if $G/H$ is torsion-free.

Now let $G$ be a $\T$-group. Recall (cf., for
example,~\cite[Chapter~6]{Segal/83}) that, by the Malcev
correspondence, there is a Lie algebra $\LieL_{ G}(\Q)$ over~$\Q$, of
dimension equal to $h( G)$, the Hirsch length of~$ G$, and an
injective mapping $\log: G\rightarrow\LieL_{ G}(\Q)$, such that ${\rm
span}_{\Q}(\log( G))=\LieL_{ G}(\Q)$. In general $L:=\log( G)$ will
fail to be a Lie subring (or even just a lattice).  However, by
replacing~$ G$ by a subgroup of finite index, if necessary, we may
assume it is a Lie subring~(\cite{GSS/88}, Section~4) and even that
$[L,L]\subseteq c!L$, where~$c$ is the nilpotency class of~$L$ (or~$
G$), i.e. that~$L$ (and $ G$) are {\sl elementarily exponentiable
(e.e.)} in Howe's nomenclature. As we are looking to study $\zirr_{
G,p}(s)$ for all but finitely many primes, there is no harm in
descending to a subgroup~$H$ of finite index in~$ G$. Indeed, for
all~$p$ and all~$n$, there is a $1-1$ -- correspondence between twist
isoclasses of irreducible characters of degree~$p^n$ and
$p$-admissible twist isoclasses of degree~$p^n$ of $\widehat{ G}_p$,
the pro-$p$-completion of~$ G$ (\cite{MartinHrushovski/04}, Lemma
8.5). However, $\widehat{ G}_p\cong\widehat{H}_p$ if $p\not|\;|G:H|$.

For a $\T$-group~$ G$, denote by~$( G^{(i)})$ the group's lower
central series, defined by~$ G^{(1)}= G$, $ G^{(i+1)}=[ G^{(i)}, G]$,
$i\geq2$. We say that~$ G\not=\{1\}$ has nilpotency class~$c$ (or is
step-$c$-nilpotent) if~$ G^{(c)}\not=\{1\}$ but~$ G^{(c+1)}=\{1\}$.
Similarly, we denote by $(L^{(i)})$ the lower central series of the
Lie algebra $L$, defined by~$L^{(1)}=L$, $L^{(i+1)}=[L^{(i)},L]$,
$i\geq2$, and we hope that there will be no confusion between group
commutators and Lie brackets. It is well-known that $\log$ induces a
bijection between $L^{(i)}$ and $ G^{(i)}$, $i\in[c]$. We write $G'$
for $G^{(2)}$, and $L'$ for $L^{(2)}$. By~$Z(L)$ we denote the centre of~$L$.

For a sublattice $M\subseteq L=\log( G)$, denote its dual
$\textup{Hom}(M,\mathbb{C}^*)$ by $\widehat{M}$ and by
$r_{M}:\widehat{L}\rightarrow\widehat{M}, \psi \mapsto \psi|_{M}$ the
restriction to~$M$. We say that~$\psi\in\widehat{L}$ is rational
on~$M$ if $r_{M}(\psi)$ is a torsion element. Clearly~$\psi$ is
rational on~$M$ if and only if~$\psi$ is rational on~$M_s$ if and only
if~$\psi(nM)=1$ for some~$n\in\N$.

Recall that the adjoint action Ad of~$ G$ on~$L$ is given by
\begin{equation}\label{adjoint formula}
l\mapsto l + [\log \gamma,l] + (\text{higher terms}),
\end{equation}
where we omitted Lie terms of degree~$\geq3$. These may be computed in
terms of the Baker-Campbell-Hausdorff-formula. This element
$F(X,Y)\in \widehat{\mathcal{L}}_{\{X,Y\}}$, the completion of the
free Lie algebra on the symbols $X$ and $Y$, provides the solution to
$$\exp(F(x,y))=\exp(x)\exp(y).$$ The Baker-Campbell-Hausdorff-formula
allows us therefore to recover the group structure on the Lie algebra
(cf.~\cite[Chapter 9]{Khukhro/98} and~\cite[V.3.4 and
IV.7]{SerreLALG/92}).  By~(\ref{adjoint formula}), the co-adjoint
action~$\Adstar$ of~$G$ on~$\widehat{L}$ is thus given by
$$\Adstar\gamma(\psi)(l)=\psi(l)\psi([\log\gamma,l])\psi(\text{higher
  terms}).$$  Given $\psi\in\widehat{L}$ we define a bi-additive,
anti-symmetric form $B_{\psi}:L\times L\rightarrow \mathbb{C}^*$ by
setting $B_{\psi}((l_1,l_2))=\psi([l_1,l_2])$. We say that a
subalgebra~$P\leq L$ {\sl polarises}~$B_{\psi}$ if~$B_{\psi}|_{P\times
  P}=1$ and it is a maximal additive subgroup with respect to this
property. Note that any such additive subgroup contains the
radical~$\Rad_{\psi}=\{l\in L|\;\psi([l,L])=1\}$ of $B_{\psi}$.

\begin{NotMyLemma}\cite[Lemmata 1--4]{Howe-nilpotent/77}\label{lemma
    howe 1234} Given $\psi\in\widehat{L}$, the isotropy subgroup
  $ G_{\psi}$ of $\psi$ under $\Adstar G$ is an e.e. subgroup of~$ G$
  and $\log  G_{\psi}=\Rad_{\psi}$. If $\psi\in\widehat{L}$ is
  $\Adstar G$-invariant, then $\psi^2$ defines a one-dimensional
  character on~$ G$. The orbit of $\psi\in\widehat{L}$ under the
  co-adjoint action~$\Adstar$ of~$ G$ is finite if and only if~$\psi$
  is rational on~$\log G'_s$.  If $\psi\in\widehat{L}$ is
  rational on~$\log G'_s$ then there are e.e. polarising
  subalgebras~$P$ for~$B_{\psi}$. They have finite index in~$L$
  satisfying~$|L:P|=|P:\Rad_{\psi}|$.
\end{NotMyLemma}

\begin{NotMyTheorem}\cite[Theorem 1(a)]{Howe-nilpotent/77}\label{theorem
    howe 1} Let~$ G$ be an e.e.~$\T$-group and set $L=\log G$,
  $L'_s=\log G'_s$. Let~$\Omega$ be a finite $\Adstar G$-orbit
  in~$\widehat{L}$, and~$\psi\in\Omega$. Let~$N$ be the period
  of~$\psi$ and assume~$N$ to be odd. A finite-dimensional irreducible
  representation~$U_{\Omega}$ may be associated to~$\Omega$ in the
  following manner: Let~$P$ be an e.e. polarising subalgebra
  for~$B_\psi$, set~$\Pi=\exp P$ and $\tilde{\psi}:=\psi|_\Pi$, a
  linear character on~$\Pi$.  Put~$U_{\Omega}:=\Ind_{\Pi}^{
  G}\tilde{\psi}$. Then the dimension of~$U_{\Omega}$
  is~$|\Omega|^{1/2}$, and the character of~$U_{\Omega}$ is
  $\xi_{\Omega}=\frac{1}{|\Omega|^{1/2}}\sum_{\phi\in\Omega}\phi$. All
  representations of the form~$\chi\otimes V$, $\chi\in\widehat{ G/
  G'_s}$, $V$ defined modulo $N G:=\exp(N\cdot L)$, $N$ odd, are
  realised in this manner.
\end{NotMyTheorem}

\begin{corollary}\label{corollary 1 howe} For almost all primes $p$,
  \begin{equation}
    \zirr_{ G,p}(s)=\sum_{\substack{\psi\in\widehat{\Ltwos}\\\psi\text{
          rational of }\\\text{$p$-power
          period}}}|L:\Rad_{\psi}|^{-s/2}|L:L_{\psi,2}|^{-1},\label{zirr howe}
\end{equation}
where
$L_{\psi,2}=\{l\in L| \psi([l,L'_s])=1\}$.
\end{corollary}

\begin{proof} For primes $p$ not dividing $2| G'_s: G'|$, Howe's
  theorem yields that we count every $p$-power degree twist-isoclass
  at least once when we sum $| G: G_\psi|^{-s/2}$ over the rational
  characters $\psi$ of ${L'_s}$ of $p$-power period. By further
  excluding finitely many primes, we may assume that $| G:
  G_\psi|=|L:\Rad_\psi|$ (cf.~\cite[Lemma 4.8]{GSS/88}, in which it is
  established that $\log$ induces an \emph{index-preserving}
  correspondence between $p$-power index subgroups of $ G$ and
  $p$-power index subalgebras of $L$ away from a finite number of
  primes). Hereby we overcount every orbit by exactly the index $| G:
  G_{\psi,2}|$, where $ G_{\psi,2}=\{\gamma\in
  G|\psi([\log\gamma,L'_s])=1 \}$ is the stabiliser of $\psi$ under
  the coadjoint action of $ G$ on the restriction of characters to
  $L'_s$. (Note that this index is always equal to~$1$ if $G$ is
  class-$2$-nilpotent.) Again at the cost of at most finitely many
  primes we may assume $| G: G_{\psi,2}|=|L:L_{\psi,2}|$.
\end{proof}

Howe's parametrisation of irreducible representations allows us to
prove\m{proposition howe}
\begin{proposition}\label{proposition howe}
  Let $ G$ be a $\T$-group. Then there are matrices $\mcS\subset\mcR$
  of homogeneous $\Q$-linear forms such that, for almost all primes
  $p$,
$$\zirr_{ G,p}(s) = P_{\mcR,\mcS,\Qp}({s/2,\dots,s/2};1,\dots,1),$$
where $P_{\mcR,\mcS,\Qp}$ is the generating function defined in
Section~\ref{subsection elementary
divisors}. Proposition~\ref{proposition elementary divisors} is
applicable for $n=h(G')$.
\end{proposition}
Theorem~\ref{theorem representations} clearly follows from this.

\begin{proof}
  We aim to express both factors in the summands of \eqref{zirr howe}
  in terms of elementary divisors of matrices of linear forms. Recall
  that $L$ is additively isomorphic to $\Z^h$, where $h$ is the Hirsch
  length of $ G$. Without loss of generality we may assume that $L'$
  is saturated in $L$, and that $L'\cap Z(L)$ is saturated in $Z(L)$
  (otherwise we disregard finitely many primes). We fix a $\Z$-basis
$$
\{x_1,\dots,x_d,\underbrace{x_{d+1},\dots,x_{d+m},\underbrace{x_{d+m+1},\dots,x_{d+n}}_{L'\cap Z(L)}}_{L'}\}
$$
of $L$ such that
\begin{align*}
\{x_{d+1},\dots,x_{d+n}\}&\text{ is a $\Z$-basis for $L'$}\\
\{x_{d+m+1},\dots,x_{d+n}\}&\text{ is a $\Z$-basis for $L'\cap Z(L)$}.
\end{align*}
The Lie bracket induces an anti-symmetric, bi-additive mapping
\begin{align*}
\beta:L/(L'\cap Z(L))\times L/(L'\cap Z(L)) &\rightarrow L'\\
\left(l_1(L'\cap Z(L)),l_2(L'\cap Z(L))\right)&\mapsto [l_1,l_2]
\end{align*}
We may express this map in terms of our chosen basis as follows. For
$1\leq i,j\leq d+m$, let
$$[x_i,x_j]=\sum_{k\in[n]}\lambda_{ijk}x_{d+k},\;
\lambda_{ijk}\in\Q.$$
Let $\mcR=(\mcR_{ij})$ denote the
$(d+m)\times(d+m)$-matrix of linear forms
$$\mcR_{ij}(\bfY)=\sum_{k\in[n]}\lambda_{ijk}Y_k\in\Q[Y_1,\dots,Y_{n}].$$
By $\mcS$ we denote the submatrix of $\mcR$ consisting of the last $m$
columns of $\mcR$. By further disregarding finitely many primes if
necessary we may assume that $p$ divides none of the denominators of
the $\lambda_{ijk}$ (so that we may reduce mod $p$), and that, given
$\bfy\in\Zp^n$, $\ol{\mcR(\bfy)}$ is not zero unless
$\ol{\bfy}=0\in\Fp^{n}$. We denote by $C_i$, $i\in[d+m]$, the
matrices of the additive maps $L/(L'\cap Z(L))\rightarrow L', x(L'\cap
Z(L))\mapsto [x,x_i]$ with respect to these bases. For a given
non-negative integer $N\in\N_0$ we identify the set $$
\Psi_{p^N}:=\left\{\psi\in\widehat{L'}|\;\text{ the period of $\psi$
equals $p^N$}\right\} $$ with
$\left(\Z^n/p^N\Z^n\right)^\times:=\primcharndprime$ by sending
$\bfl=(l_1,\dots,l_{n})\in\left(\Z^n/p^N\Z^n\right)^\times$ to
$$\psi(b_1,\dots,b_{n})=\exp\left(\frac{2\pi i\sum_{i\in[n]}
    l_ib_i}{p^N}\right),\; \bfb=(b_1,\dots,b_{n})\in\widehat{L'}.$$
With these identifications we obtain
\begin{align}
  \bfk=(k_1,\dots,k_{d+m})\in\Rad_\psi/Z(L)\Leftrightarrow& \forall
  \gamma\in G:\;\psi\left(\left[\sum
  k_ix_i,\log\gamma\right]\right)=1\nonumber\\ \Leftrightarrow&\forall j\in[d+m]:\;\bfk\;
  C_j\bfl^t\equiv 0 \mod p^N\nonumber\\
  \Leftrightarrow&\bfk\; \mcR(\bfl)\equiv 0 \mod p^N\label{congruence
  R}
\end{align}
and
\begin{align}
  \bfk=(k_1,\dots,k_{d+m})\in L_{\psi,2}/Z(L)\Leftrightarrow& \forall
  \gamma\in G':\;\psi\left(\left[\sum
  k_ix_i,\log\gamma\right]\right)=1\nonumber\\ \Leftrightarrow&\forall j\in[d+1,d+m]:\;\bfk\;
  C_j\bfl^t\equiv 0 \mod p^N\nonumber\\
  \Leftrightarrow&\bfk\; \ \mcS(\bfl)\equiv 0 \mod
  p^N\label{congruence S}.
\end{align}
In order to use the congruence conditions \eqref{congruence R} and
\eqref{congruence S} for an effective computation
of~$\zirr_{ G,p}(s)$, we need to enumerate the ($p$-parts of) the
elementary divisors of the matrices~$\mcR(\bfl)$ and $\mcS(\bfl)$
as~$\bfl$ runs through the sets~$\Psi_{p^N}$, $N\in\N$.  By
Corollary~\ref{corollary 1 howe} we may write
\begin{align*}
  \zirr_{G,p}(s)=&\sum_{\substack{N\in\N_0\\\psi\in\Psi_{p^N}}}|L:\Rad_\psi|^{-s/2}|L:L_{\psi,2}|^{-1}\\
  =&\sum_{\substack{N\in\N_0\\\bfm\in\N_0^{d+m},\;\bfn\in\N_0^m}}
  \mcN_{N,\bfm,\bfn}q^{-\sum_{i\in[d+m]}(N-m_i)^{s/2}-\sum_{j\in[m]}(N-n_j)}\\
  =&P_{\mcR,\mcS,\Qp}(s/2,\dots,s/2;1,\dots,1).
\end{align*}
\end{proof}

\begin{acknowledgements}  
  I started work on this paper when I was a guest at the
  Max-Planck-Institut f\u r Mathematik in Bonn, Germany. I continued
  it as a guest at the Chaire de Math\'ematiques Alg\'ebriques
  Discr\`etes at the EPF Lausanne, Switzerland. I gratefully
  acknowledge both institution's support and hospitality. I am
  particularly obliged to Laurent Bartholdi, my host in Lausanne, and
  to the Fond National Suisse, which jointly funded my stay there.
  Thanks are due to Benjamin {Klopsch} for frequent inspiring
  discussions. I wish to thank Benjamin Martin for making Hrushovski's
  and his preprint~\cite{MartinHrushovski/04} available to me, and to
  Willem Veys for alerting me to \cite{VeysZG/06}.  Mark Berman's and
  Ulrich Voll's comments on previous versions helped a great deal to
  improve the presentation of this paper. I am particularly grateful
  to Mark Berman for his help with proof-reading.
\end{acknowledgements}
\providecommand{\bysame}{\leavevmode\hbox to3em{\hrulefill}\thinspace}
\providecommand{\MR}{\relax\ifhmode\unskip\space\fi MR }
\providecommand{\MRhref}[2]{%
  \href{http://www.ams.org/mathscinet-getitem?mr=#1}{#2}
}
\providecommand{\href}[2]{#2}

\end{document}